\documentclass[11pt]{amsart} 
\usepackage{amsmath,amssymb}
\usepackage[dvips]{graphicx}  
\usepackage{psfrag}
\def\COMMENT#1{}
\def\TASK#1{}

\def\noproof{{\unskip\nobreak\hfill\penalty50\hskip2em\hbox{}\nobreak\hfill%
        $\square$\parfillskip=0pt\finalhyphendemerits=0\par}\goodbreak}
\def\endproof{\noproof\bigskip}
\newdimen\margin   
\def\textno#1&#2\par{%
    \margin=\hsize
    \advance\margin by -4\parindent
           \setbox1=\hbox{\sl#1}%
    \ifdim\wd1 < \margin
       $$\box1\eqno#2$$%
    \else
       \bigbreak
       \hbox to \hsize{\indent$\vcenter{\advance\hsize by -3\parindent
       \sl\noindent#1}\hfil#2$}%
       \bigbreak
    \fi}
\def\proof{\removelastskip\penalty55\medskip\noindent{\bf Proof. }}

\def\eps{\varepsilon}

\newcommand{\msc}[1]{\begin{center}MSC2000: #1.\end{center}}
\newtheorem{firstthm}{Proposition}
\newtheorem{thm}[firstthm]{Theorem}
\newtheorem{prop}[firstthm]{Proposition}
\newtheorem{lemma}[firstthm]{Lemma}

\newtheorem{conj}[firstthm]{Conjecture}
\newtheorem{claim}[firstthm]{Claim}

\newtheorem{question}[firstthm]{Question}

\begin{document}
\title{Hamilton decompositions of regular tournaments}
\author{Daniela K\"uhn, Deryk Osthus and Andrew Treglown}
\thanks {The authors were supported by the EPSRC, grant no.~EP/F008406/1.}
\date{\today} 

\begin{abstract}
We show that every sufficiently large regular tournament can almost completely be decomposed into edge-disjoint 
Hamilton cycles. More precisely, for each $\eta >0$ every regular tournament $G$ of sufficiently large
order $n$ contains at least $(1/2-\eta)n$ edge-disjoint Hamilton cycles. This gives an approximate solution
to a conjecture of Kelly from 1968. Our result also extends to almost regular tournaments.
\end{abstract}

\maketitle
\msc{5C20, 5C35, 5C45}
\section{Introduction}\label{sec1}
A Hamilton decomposition of a graph or digraph $G$ is a set of edge-disjoint Hamilton cycles which together
cover all the edges of~$G$.
The topic has a long history but some of the main questions remain open. 
In 1892, Walecki  showed that the edge set of the complete graph $K_n$ on $n$ vertices has a
Hamilton decomposition if $n$ is odd (see e.g.~\cite{abs,lucas} for the construction). 
If $n$ is even, then $n$ is not a factor of ${n \choose 2}$, so clearly $K_n$ does not 
have such a decomposition.
Walecki's result implies that a complete digraph $G$ on $n$ vertices has a Hamilton decomposition if $n$ is odd. 
More generally, Tillson~\cite{till} proved 
that a complete digraph $G$ on $n$ vertices has a Hamilton decomposition if and only if $n \not = 4,6$.  
 
A tournament is an orientation of a complete graph. We say that a tournament is {\it regular} if every vertex has 
equal in- and outdegree. Thus
regular tournaments contain an odd number $n$ of vertices and each vertex has in- and outdegree  $(n-1)/2$. 
The following beautiful conjecture of Kelly~(see e.g.~\cite{bang,bondy,moon}), which has attracted much attention,
states that every regular tournament has a Hamilton decomposition:
\begin{conj}[Kelly] \label{kelly}
Every regular tournament on $n$ vertices can be decomposed
into $(n-1)/2$ edge-disjoint Hamilton cycles.
\end{conj}
In this paper we prove an approximate version of Kelly's conjecture. 
\begin{thm}\label{main1} For every $\eta >0$ there exists an integer $n_0$ 
so that every regular tournament on $n \geq n_0$ vertices  contains at least $(1/2-\eta)n$ edge-disjoint Hamilton
cycles.
\end{thm}
In fact, we prove the following stronger result, where
we consider orientations of almost complete graphs which are almost regular.
An \emph{oriented graph} is obtained from an undirected graph by orienting its edges. So it has at most one edge
between every pair of vertices, whereas a digraph may have an edge in each direction.
\begin{thm}\label{main}
For every $\eta_1 >0$ there exist $n_0= n_0 (\eta_1)$ and $\eta _2=\eta _2 (\eta_1) >0$ 
such that the following holds. Suppose that $G$ is an oriented graph on $n\geq n_0$ vertices such that
every vertex in $G$ has in- and outdegree at least $(1/2-\eta _2) n$.
Then $G$ contains at least $(1/2-\eta_1)n$ edge-disjoint Hamilton cycles. 
\end{thm}
The \emph{minimum semidegree} $\delta^0 (G)$ of an oriented graph 
$G$ is the minimum of its minimum outdegree and its 
minimum indegree. So the minimum semidegree of a regular tournament on $n$ vertices is $(n-1)/2$.
Most of the previous partial results towards Kelly's conjecture have been obtained by giving bounds on the  
minimum semidegree of an oriented graph which guarantees a Hamilton cycle.
This approach was first used by Jackson~\cite{jackson}, who showed that
every regular tournament on at least 5 vertices contains a Hamilton cycle and a Hamilton path which are edge-disjoint.
Zhang~\cite{zhang} then showed that every such tournament contains two edge-disjoint Hamilton cycles.
Improved bounds on the value of $\delta^0 (G)$ which forces a Hamilton cycle were then found by
Thomassen~\cite{tom1}, H\"aggkvist~\cite{HaggkvistHamilton},
H\"aggkvist and Thomason~\cite{haggtom} as well as Kelly, K\"uhn and Osthus~\cite{kelly}.
Finally, Keevash, K\"uhn and Osthus~\cite{kko} showed that every sufficiently large oriented graph~$G$ on $n$ vertices
with $\delta ^0 (G) \geq (3n-4)/8$ contains a Hamilton cycle. This bound on $\delta^0(G)$ is best possible
and confirmed a conjecture of 
H\"aggkvist~\cite{HaggkvistHamilton}. Note that this result implies that every sufficiently large 
regular tournament on $n$ vertices contains at least $n/8$ edge-disjoint Hamilton cycles. 
This was the best bound so far towards Kelly's conjecture. 

Kelly's conjecture has also been verified for $n \le 9$ by Alspach (see the survey~\cite{bt}).
A result of Frieze and Krivelevich~\cite{fk} states that Theorem~\ref{main} holds for
`quasi-random' tournaments. As indicated below, we will build on some of their
ideas in the proof of Theorem~\ref{main}.

It turns out that Theorem~\ref{main} can be generalized even further:
any large almost regular oriented graph on $n$ vertices whose in- and outdegrees are all a little larger than 
$3n/8$ can almost be decomposed into Hamilton cycles.
The corresponding modifications to the proof of Theorem~\ref{main} are described in Section~\ref{38}.
We also discuss some further open questions in that section.

Jackson~\cite{jackson} also introduced the following bipartite version of Kelly's conjecture
(both versions are also discussed e.g.~in the Handbook article by Bondy~\cite{bondy}).
A \emph{bipartite tournament} is an orientation of a complete bipartite graph.
\begin{conj}[Jackson] \label{kellybip}
Every regular bipartite tournament has a Hamilton decomposition.
\end{conj}
An undirected version of Conjecture~\ref{kellybip} was proved
independently by Auerbach and Laskar~\cite{auerbach}, as well as Hetyei~\cite{hetyei}.
However, a bipartite version of Theorem~\ref{main} does not hold, because there are almost regular bipartite
tournaments which do not even contain a single Hamilton cycle. 
(Consider for instance the following `blow-up' of a 4-cycle: the vertices are split into 4 parts $A_0,\dots,A_3$
whose sizes are almost but not exactly equal, and we have all edges from $A_i$ to $A_{i+1}$, with indices modulo 4.)

Kelly's conjecture has been generalized in several directions. For instance, given 
an oriented graph $G$, define its \emph{excess} by
$$
{\rm ex}(G):=\sum_{v \in V(G)} \max \{ d^+(v)-d^-(v),0 \},
$$ 
where $d^+(v)$ denotes the number of outneighbours of the vertex $v$, and $d^-(v)$ the number of
its inneighbours. Pullman (see e.g.~Conjecture~8.25 in~\cite{bondy}) conjectured that
if $G$ is an oriented graph such that $d^+(v)+d^-(v)=d$ for all vertices $v$ of $G$, where $d$ is
odd, then $G$ has a decomposition into ${\rm ex}(G)$
directed paths. To see that  this would imply Kelly's conjecture, let $G$ be the oriented
graph obtained from a regular tournament by deleting a vertex.
Another generalization was made by
Bang-Jensen and Yeo~\cite{bangyeo}, who conjectured that every $k$-edge-connected tournament 
has a decomposition into $k$ spanning strong digraphs.
 
In~\cite{tom1}, Thomassen also formulated the following weakening of Kelly's conjecture.
\begin{conj}[Thomassen] \label{thomconj}
If $G$ is a regular tournament on $2k+1$ vertices and $A$ is any set of at most
$k-1$ edges of $G$, then $G-A$ has a Hamilton cycle.
\end{conj}
In~\cite{chvatal}, we proved a result on the existence of
Hamilton cycles in `robust expander digraphs' which implies Conjecture~\ref{thomconj} 
for large tournaments (see~\cite{chvatal} for details).
\cite{tom1} also contains the related conjecture
that for any $\ell \ge 2$, there is an $f(\ell)$ so that every strongly $f(\ell)$-connected 
tournament contains $\ell$ edge-disjoint Hamilton cycles.

Further support for Kelly's conjecture was also provided by Thomassen~\cite{thom2}, who showed that 
the edges of every regular tournament on $n$ vertices can be covered by $12n$ Hamilton cycles.
In~\cite{cyclesurvey} the first two authors observed that one can use Theorem~\ref{main} to reduce this to $(1/2+o(1))n$ Hamilton cycles.
A discussion of further recent results about Hamilton cycles in directed graphs can be found in the survey~\cite{cyclesurvey}.

It seems likely that the techniques developed in this paper will also be useful in solving further problems.
In fact, Christofides, K\"uhn and Osthus~\cite{cko} used similar ideas to prove  approximate versions 
of the following two long-standing conjectures of Nash-Williams~\cite{nash1, nash2}:%
\begin{conj}[Nash-Williams~\cite{nash1}]
Let $G$ be a $2d$-regular graph on at most $4d+1$ vertices, where $d \ge 1$. 
Then $G$ has a Hamilton decomposition. 
\end{conj}
\begin{conj}[Nash-Williams~\cite{nash2}] \label{nwconj2}
Let $G$ be a graph on $n$ vertices with minimum degree at least $n/2$. 
Then $G$ contains $n/8+o(n)$ edge-disjoint Hamilton cycles.
\end{conj}
(Actually, Nash-Williams initially formulated Conjecture~\ref{nwconj2} with the term $n/8$ replaced
by $n/4$, but Babai found a counterexample to this.)

Another related problem was raised by Erd\H{o}s (see~\cite{tom1}), who asked whether almost all tournaments $G$
have at least $\delta^0(G)$ edge-disjoint Hamilton cycles. Note that an affirmative answer would not directly imply
that Kelly's conjecture holds for almost all regular tournaments, which would of course be an interesting result in itself.
There are also a number of corresponding questions for random undirected graphs (see e.g.~\cite{fk}).

After giving an outline of the argument in the next section, we will state
a directed version of the Regularity lemma and some related results in Section~\ref{3}.
Section~\ref{4} contains statements and proofs of several auxiliary results, mostly 
on (almost) $1$-factors in (almost) regular oriented graphs.
The proof of Theorem~\ref{main} is given in Section~\ref{5}.
A generalization of Theorem~\ref{main} to oriented graphs with smaller degrees is
discussed in Section~\ref{38}.

\section{Sketch of the proof of Theorem~\ref{main}} \label{sketch}

Suppose we are given a regular tournament $G$ on $n$ vertices and our aim is to `almost' decompose it into Hamilton cycles.
One possible approach might be the following: first remove a spanning regular oriented subgraph $H$ whose degree $\gamma n$ 
satisfies $\gamma \ll 1$. Let $G'$ be the remaining oriented subgraph of $G$.
Now consider a decomposition of $G'$ into $1$-factors $F_1,\dots,F_r$ (which clearly exists).
Next, try to transform each $F_i$ into a Hamilton cycle by removing some of its edges and adding some suitable edges of $H$.
This is of course impossible if many of the $F_i$ consist of many cycles.
However, an auxiliary result of Frieze and Krivelevich in~\cite{fk} implies that we can `almost' decompose $G'$ so that
each $1$-factor $F_i$ consists of only a few cycles.

If $H$ were a `quasi-random' oriented graph, then (as in~\cite{fk}) one could use it to successively `merge' the cycles
of each $F_i$ into Hamilton cycles using a `rotation-extension' argument: 
delete an edge of a cycle $C$ of $F_i$ to obtain a path $P$ from $a$ to $b$, say.
If there is an edge of $H$ from $b$ to another cycle $C'$ of $F_i$, then extend $P$ to include the vertices of $C'$
(and similarly for $a$). Continue until there is no such edge. Then (in $H$) the current endvertices of the path $P$ have many
neighbours on $P$. One can use this together with the quasi-randomness of $H$ to transform $P$ into a cycle with the same vertices 
as $P$. Now repeat this, until we have merged all the cycles into a single (Hamilton) cycle.
Of course, one has to be careful to maintain the quasi-randomness of $H$ in carrying out this `rotation-extension' process
for the successive $F_i$ (the fact that $F_i$ contains only few cycles is important for this). 

The main problem is that $G$ need not contain such a spanning `quasi-random' subgraph $H$.
So instead, in Section~\ref{applyDRL} we use Szemer\'edi's regularity lemma to decompose $G$ into quasi-random subgraphs.
We then choose both our $1$-factors $F_i$ and the graph $H$ according to the structure of this decomposition.
More precisely, we apply a directed version of Szemer\'edi's regularity lemma to obtain a partition of the vertices of
$G$ into a bounded number of clusters $V_i$ so that almost all of the bipartite subgraphs spanned by ordered pairs of clusters
are quasi-random (see Section~\ref{3.3} for the precise statement). 
This then yields a reduced digraph $R$, whose vertices correspond to the clusters, with an 
edge from one cluster $U$ to another cluster $W$ if the edges from $U$ to $W$ in $G$ form a quasi-random graph.
(Note that $R$ need not be oriented.) We view $R$ as a weighted digraph whose edge weights are the densities of the 
corresponding ordered pair of clusters. We then obtain an unweighted multidigraph $R_m$ from $R$ as follows:
given an edge $e$ of $R$ joining a cluster $U$ to $W$, replace it with $K=K(e)$ copies of $e$, 
where $K$ is approximately proportional to the density of the ordered pair $(U,W)$.
It is not hard to show that $R_m$ is approximately regular (see Lemma~\ref{multimin}).
If $R_m$ were regular, then it would have a decomposition into $1$-factors,
but this assumption may not be true. However, we can show that $R_m$ can `almost' be decomposed into `almost' $1$-factors.
In other words, there exist edge-disjoint collections $\mathcal F_1, \dots , \mathcal F_r$ of vertex-disjoint
cycles in $R_m$ such that each $\mathcal F_i$ covers almost all of the clusters in $R_m$ (see Lemma~\ref{multifactor1}).

Now we choose edge-disjoint oriented spanning subgraphs $C_1,\dots,C_r$ of $G$ so that each $C_i$ corresponds to $\mathcal F_i$.
For this, consider an edge $e$ of $R$ from $U$ to $W$ and suppose 
for example that $\mathcal F_1$, $\mathcal F_2$ and $\mathcal F_8$ are the only $\mathcal F_i$ containing copies 
of $e$ in $R_m$. 
Then for each edge of $G$ from $U$ to $W$ in turn, we assign it to one of $C_1$, $C_2$ and $C_8$ with equal probability.
Then with high probability, each $C_i$ consists of bipartite quasi-random oriented graphs which together
form a disjoint union of `blown-up' cycles. Moreover, we can arrange that all the vertices have degree close to $\beta m$
(here $m$ is the cluster size and $\beta$ a small parameter which does not depend on $i$).
We now remove a small proportion of the edges from $G$ (and thus from each $C_i$) to form oriented subgraphs 
$H_1^+,H_1^-,H_2,H_{3,i},H_4,H_{5,i}$ of $G$, where $1 \le i \le r$.
Ideally, we would like to show that each $C_i$ can almost be decomposed into Hamilton cycles.
Since the $C_i$ are edge-disjoint, this would yield the required result.

One obvious obstacle is that the $C_i$ need not be spanning subgraphs of $G$
(because of the exceptional set $V_0$ returned by the regularity lemma 
and because the $\mathcal F_i$ are not spanning.)
So in Section~\ref{sec:incorp} we add suitable edges between $C_i$ and the leftover vertices to form edge-disjoint oriented
spanning subgraphs $G_i$ of $G$ where every vertex has degree close to $\beta m$. 
(The edges of $H_1^-$ and $H_1^+$ are used in this step.)
But the distribution of the edges added in this step may be somewhat `unbalanced',
with some vertices of $C_i$ sending out or receiving too many of them.
In fact, as discussed at the beginning of Section~\ref{skel},
we cannot even guarantee that $G_i$ has a single $1$-factor. We overcome this new difficulty by 
adding carefully chosen further edges (from $H_2$ this time) to each $G_i$ which compensate the above
imbalances. 

Once these edges have been added, 
in Section~\ref{nicefactor} we can use the max-flow min-cut theorem
to almost decompose each $G_i$ into $1$-factors $F_{i,j}$. 
(This is one of the points where we use the fact that the $C_i$ consist of quasi-random graphs which form
a union of blown-up cycles.) Moreover, (i) the number of cycles in each 
of these $1$-factors is not too large and (ii) most of the cycles inherit the structure of
$\mathcal F_i$. More precisely, (ii) means that most vertices $u$ of $C_i$ have the following property:
let $U$ be the cluster containing $u$ and let $U^+$ be the successor of $U$ in $\mathcal F_i$. Then 
the successor $u^+$ of $u$ in $F_{i,j}$ lies in $U^+$.

In Section~\ref{4.6} we can use (i) and (ii) to merge the cycles of each $F_{i,j}$ into a 
$1$-factor $F'_{i,j}$ consisting only of a bounded number of cycles
-- for each cycle $\mathcal C$ of $\mathcal F_i$, all the vertices 
of $G_i$ which lie in clusters of $\mathcal C$ will lie in the same cycle of $F_{i,j}'$.
We will apply a rotation-extension argument for this, where the additional edges (i.e.~those not in $F_{i,j}$)
come from $H_{3,i}$.
Finally, in Section~\ref{merging} we will use the fact that $R_m$ contains many short paths to merge each $F'_{i,j}$ into a single 
Hamilton cycle. The additional edges will come from $H_4$ and $H_{5,i}$ this time.

\section{Notation and the Diregularity lemma} \label{3}
\subsection{Notation}
Throughout this paper we omit floors and ceilings whenever this does not affect the argument.
Given a graph $G$, we denote the degree of a vertex $x \in V(G)$ by $d_G (x)$ and the maximum degree of $G$ by
$\Delta (G)$. Given two vertices $x$ and $y$ of a digraph $G$, we write $xy$ for
the edge directed from $x$ to $y$. We denote by $N^+ _G (x)$ the set of all outneighbours of~$x$.
So $N^+ _G (x)$ consists of all those $y \in V(G)$ for which
$xy \in E(G)$. We have an analogous definition for $N^-_G (x).$ Given a multidigraph $G$, we denote
by $N^+ _{G} (x)$ the {\emph{multiset}} of vertices where a vertex $y \in V(G)$ appears $k$ times in 
$N_G ^+ (x)$ if $G$ contains precisely $k$ edges from~$x$ to~$y$. Again, we have an analogous
definition for $N^-_G (x)$. We will write $N^+ (x)$ for example, if this is unambiguous.
Given a vertex $x$ of a digraph or multidigraph $G$, we write $d^+ _G (x):=|N^+(x)|$ for the
outdegree of $x$, $d^- _G(x):=|N^-(x)|$ for its indegree and $d(x):=d^+(x)+d^-(x)$ for
its degree. The maximum of the maximum outdegree
$\Delta ^+ (G)$ and the maximum indegree $\Delta ^- (G) $ is denoted by $\Delta ^0 (G)$.
The  \emph{minimum semidegree} $\delta^0 (G)$ of $G$ is the minimum of its minimum outdegree
$\delta ^+ (G)$ and its minimum indegree $\delta ^- (G)$. Throughout the paper we will use
$d^{\pm} _G (x)$, $\delta ^{\pm} (G)$ and $N^{\pm} _G (x)$ as `shorthand' notation. For example, 
$\delta ^{\pm} (G) \geq \delta ^{\pm} (H)/2$ is read as 
$\delta ^+ (G) \geq \delta ^+ (H)/2$ and
$\delta ^- (G) \geq \delta ^- (H)/2$.

A \emph{1-factor} of a multidigraph~$G$ is a collection of vertex-disjoint cycles in~$G$ which
together cover all the vertices of~$G$.
Given $A,B \subseteq V(G)$, we write $e_G (A,B)$ to denote the number of edges in $G$ with startpoint in $A$
and endpoint in $B$. Similarly, if $G$ is an undirected graph, we write $e_G (A,B)$ for the number of
all edges between $A$ and~$B$. Given a multiset $X$ and a set $Y$ we define $X \cap Y$ to be the multiset
where $x$  appears as an element precisely $k$ times in $X \cap Y$ if $x \in X$, $x \in Y$ and $x$
appears precisely $k$ times in $X$. We write $a=b\pm \eps$ for $a\in [b-\eps,b+\eps]$.


\subsection{A Chernoff bound}
We will often use the following Chernoff bound for binomial and hypergeometric
distributions (see e.g.~\cite[Corollary 2.3 and Theorem 2.10]{Janson&Luczak&Rucinski00}).
Recall that the binomial random variable with parameters $(n,p)$ is the sum
of $n$ independent Bernoulli variables, each taking value $1$ with probability $p$
or $0$ with probability $1-p$.
The hypergeometric random variable $X$ with parameters $(n,m,k)$ is
defined as follows. We let $N$ be a set of size $n$, fix $S \subseteq N$ of size
$|S|=m$, pick a uniformly random $T \subseteq N$ of size $|T|=k$,
then define $X=|T \cap S|$. Note that $\mathbb{E}X = km/n$.

\begin{prop}\label{chernoff}
Suppose $X$ has binomial or hypergeometric distribution and $0<a<3/2$. Then
$\mathbb{P}(|X - \mathbb{E}X| \ge a\mathbb{E}X) \le 2 e^{-\frac{a^2}{3}\mathbb{E}X}$.
\end{prop}


\subsection{The Diregularity lemma} \label{3.3}
In the proof of Theorem~\ref{main} we will use the directed version of Szemer\'edi's Regularity lemma.
Before we can state it we need some more notation and definitions.
The \emph{density} of an undirected bipartite graph $G$ with vertex classes~$A$ and~$B$ is
defined to be
$$d_G(A,B):=\frac{e_G(A,B)}{|A||B|}.$$
We will write $d(A,B)$ if this is unambiguous.
Given any $\eps, \eps '>0$, we say that $G$ is  \emph{$[\eps,\eps']$-regular} if  for all sets
$X \subseteq A$ and $Y \subseteq B$ with $|X|\ge \eps |A|$ and
$|Y|\ge \eps |B|$ we have $|d(A,B)-d(X,Y)|< \eps '$. 
In the case when $\eps =\eps '$ we say that $G$ is \emph{$\eps$-regular}.

Given $d \in [0,1)$ we say that $G$ is \emph{$(\eps,d)$-super-regular} if all sets
$X \subseteq A$ and $Y \subseteq B$ with $|X|\ge \eps |A|$ and $|Y|\ge \eps |B|$
satisfy $d(X,Y)= d\pm \eps$ and, furthermore, if $ d_G (a)= (d\pm \eps)|B|$ for
all $a \in A$ and $ d_G (b)=(d\pm \eps )|A|$ for all $b \in B$. Note that this is a slight
variation of the standard definition.

Given disjoint vertex sets~$A$ and~$B$ in a digraph~$G$, we write $(A,B)_G$ for the
oriented bipartite subgraph of~$G$
whose vertex classes are~$A$ and~$B$ and whose edges are all the edges from~$A$ to~$B$ in~$G$.
We say $(A,B)_G$ is \emph{$[\varepsilon , \eps ']$-regular and has density~$d'$} if this holds
for the underlying undirected bipartite graph of $(A,B)_G$.
(Note that the ordering of the pair $(A,B)_G$ is important here.)
In the case when $\eps = \eps '$ we say that \emph{$(A,B)_G$ is $\eps$-regular and has density $d'$}.
Similarly, given $d \in [0,1)$ we say $(A,B)_G$ is \emph{$(\eps ,d) $-super-regular} if this holds 
for the underlying undirected bipartite graph.
 
The Diregularity lemma is a variant of the Regularity lemma for digraphs due to Alon
and Shapira~\cite{alon}. Its proof is similar to the undirected version.
We will use the degree form of the Diregularity lemma which can be derived from the standard
version in the same manner as the undirected degree form (see~\cite{survey} for a sketch of
the latter).

\begin{lemma}[Degree form of the Diregularity lemma]\label{dilemma}
For every $\varepsilon\in (0,1)$ and every integer~$M'$ there are integers~$M$ and~$n_0$
such that if~$G$ is a digraph on $n\ge n_0$ vertices and
$d\in[0,1]$ is any real number, then there is a partition of the vertex set of~$G$ into
$V_0,V_1,\ldots,V_L$ and a spanning subdigraph~$G'$ of~$G$ such that the following holds:
\begin{itemize}
\item $M'\le L\leq M$,
\item $|V_0|\leq \varepsilon n$,
\item $|V_1|=\dots=|V_L|=:m$,
\item $d^\pm_{G'}(x)>d^\pm_G(x)-(d+\varepsilon)n$ for all vertices $x\in V(G)$,
\item for all $i=1,\dots,L$ the digraph $G'[V_i]$ is empty,	
\item for all $1\leq i,j\leq L$ with $i\neq j$ the pair $(V_i,V_j)_{G'}$ is $\varepsilon$-regular
and has density either~$0$ or at least~$d$.
\end{itemize}
\end{lemma}
We call $V_1, \dots, V_L$ \emph{clusters}, $V_0$ the \emph{exceptional set} and the
vertices in~$V_0$ \emph{exceptional vertices}. We refer to~$G'$ as the \emph{pure digraph}.
The last condition of the lemma says that all pairs of clusters are $\eps$-regular in both
directions (but possibly with different densities).
The {\it reduced digraph~$R$ of~$G$ with parameters $\varepsilon$, $d$ and~$M'$} is the digraph whose 
vertices are $V_1, \dots , V_L$ and in which $V_i V_j$ is an edge precisely when $(V_i,V_j)_{G'}$
is $\varepsilon$-regular and has density at least~$d$.

The next result shows that we can partition the set of edges of an $\eps$-(super)-regular pair into edge-disjoint
subgraphs such that each of them is still (super)-regular.

\begin{lemma}\label{split}
Let $0< \eps  \ll d_0 \ll 1$ and suppose $K\ge 1$. Then there exists an integer
$m_0=m_0 (\eps,d_0, K)$ such that for all $d\ge d_0$ the following holds.
\begin{itemize}
\item[{\rm (i)}] Suppose that
$G=(A,B)$ is an $\eps$-regular pair of density $d$ where $|A|=|B|=m\geq m_0$.
Then there are $\lfloor K\rfloor$ edge-disjoint spanning subgraphs $S_1, \dots, S_{\lfloor K\rfloor}$ of
$G$ such that each $S_i$ is $[\eps, 4\eps /K]$-regular of density $(d\pm 2\eps)/K$.
\item[{\rm (ii)}] If $K=2$ and $G=(A,B)$ is $(\eps,d)$-super-regular with $|A|=|B|=m\geq m_0$.
then there are two edge-disjoint spanning subgraphs $S_1$ and $S_2$ of $G$ such that
each $S_i$ is $(2\eps,d/2)$-super-regular.
\end{itemize}
\end{lemma}
\proof We first prove~(i). Suppose we have chosen $m_0$ sufficiently large. 
Initially set $E(S_i)= \emptyset$ for each $i=1, \dots , \lfloor K \rfloor$. We consider each edge of $G$ in turn and
add it to each $E(S_i)$ with probability $1/K$, independently of all other edges of~$G$. 
So the probability that $xy$ is added to none of the $S_i$ is $1-\lfloor K\rfloor/K$.
Moreover, $\mathbb E (e(S_i)) =  e(G)/K=d m^2/K$. 

Given $X \subseteq A$ and $Y \subseteq B$ with $|X|,|Y|\geq \eps m$ we have that $|d_G(X,Y)-d|<\eps$.
Thus
$$ \frac{1}{K}(d- \eps)|X||Y|<\mathbb E(e_{S_i} (X,Y)) < \frac{1}{K}(d+ \eps)|X||Y|$$
for each $i$. Proposition~\ref{chernoff} for the binomial distribution implies that with high probability
$(d - 2\eps)|X||Y|/K < e_{S_i} (X,Y)<(d + 2\eps)|X||Y| /K$ for each $i \leq \lfloor K \rfloor$ and
every $X \subseteq A$ and $Y \subseteq B$ with $|X|,|Y|\geq \eps m$.
Such~$S_i$ are as required in~(i).\COMMENT{We are looking at roughly $(2^m)^2=4^m$ pairs of sets $X,Y$ here (for each $i=1, \dots , \lfloor
K \rfloor $). So in total using Chernoff $\leq K 4^m $ times. For each $Z:=e_{S_i} (X,Y)$, we have that 
$\mathbb P (|Z-\mathbb E Z|\geq \eps \mathbb E Z)\leq 2e^{-\frac{\eps ^2}{3} \mathbb E Z} \leq 2 e^{-cm^2}$ for some $c>0$ (since $\mathbb 
E Z\geq (d-\eps)\eps ^2 m^2/K$). 
Now $2K4^m e^{-cm^2} \ll 1$ as $m$ sufficiently large, as desired.}

The proof of~(ii) is similar. Indeed, as in (i) one can show that with high probability any 
$X\subseteq A$ and $Y \subseteq B$ with $|X|,|Y|\geq \eps m$ satisfy
$d_{S_i} (X,Y)=d/2 \pm 2 \eps$ (for $i=1,2$).
Moreover, each 
vertex $a\in A$ satisfies $\mathbb E(d_{S_i}(a))=d_G(a)/2=(d\pm \eps)m/2$ (for $i=1,2$) and similarly
for the vertices in~$B$. So again Proposition~\ref{chernoff} for the binomial distribution implies that
with high probability $d_{S_i}(a)=(d/2\pm 2\eps)m$ for all $a\in A$ and $d_{S_i}(b)=(d/2\pm 2\eps)m$ for all $b\in B$.
Altogether this shows that with high probability both $S_1$ and $S_2$ are $(2\eps,d/2)$-super-regular.\COMMENT{Same as in case (i). We have 
to consider the $2m$ degrees as well. Now $\mathbb P (|d_{S_i} (x)-\mathbb E d_{S_i} (x)|\geq \eps\mathbb E  d_{S_i} (x) ) \leq 2 e ^{-c'm} $
for some $c'>0$. But as $4^{m+1} e^{-cm^2} +4m e^{-c'm} \ll 1$ we are fine.}
\endproof 

Suppose $0<1/M'\ll \eps \ll \beta \ll d \ll 1$ and let $G$ be a digraph. Let $R$ and $G'$ denote the
reduced digraph and pure digraph respectively, obtained 
by applying Lemma~\ref{dilemma} to $G$ with parameters $\eps, d$ and $M'$. For each edge $V_iV_j$ of~$R$
we write $d_{i,j}$ for the density of $(V_i, V_j)_{G'}$. (So $d_{i,j}\geq d$.) The
\emph{reduced multidigraph} $R_m$ of $G$ with parameters 
$\eps, \beta, d$ and $M'$ is obtained from $R$ 
by setting $V(R_m):=V(R)$ and adding $\lfloor d_{i,j} /\beta \rfloor$ directed
edges from $V_i$ to $V_j$ whenever $V_iV_j\in E(R)$.

We will always consider the reduced multidigraph $R_m$ of a digraph $G$ whose order is
sufficiently large in order to apply Lemma~\ref{split} to any pair $(V_i, V_j)_{G'}$ of clusters
with $V_iV_j\in E(R)$. Let $K:= d_{i,j}/\beta$ and $S_{i,j,1}, \dots, S_{i,j,\lfloor K\rfloor}$
be the spanning subgraphs of $(V_i,V_j)_{G'}$ obtained from Lemma~\ref{split}. 
(So each $S_{i,j,k}$ is $\eps$-regular of density $\beta\pm \eps$.)
Let $(V_iV_j)_1, \dots , (V_iV_j)_{\lfloor K\rfloor}$ denote the directed edges from $V_i$ to $V_j$ in $R_m$. 
We associate each $(V_iV_j)_k$ with the edges in $S_{i,j,k}$.  

\begin{lemma}\label{multimin}
Let $0 < 1/M' \ll \eps \ll \beta \ll d \ll c_1 \leq c_2 <1$ and let $G$ be a digraph of
sufficiently large order $n$ with $\delta ^0 (G) \geq c_1n$ and $\Delta ^0 (G) \leq c_2n$. 
Apply Lemma~\ref{dilemma} with parameters $\eps, d$ and $M'$ to obtain a pure digraph $G'$
and a reduced digraph $R$ of $G$. Let $R_m$ denote the reduced multidigraph of $G$ with
parameters $\eps, \beta, d$ and $M'$. Then
$$\delta ^0 (R_m) > (c_1-3d)\frac{|R_m|}{\beta} \text{ and } \Delta ^0 (R_m)< (c_2 +2 \eps)\frac{|R_m|}{\beta}.$$
\end{lemma}
\noindent
Note the corresponding upper bound would not hold if we considered $R$ instead of $R_m$ here.
\proof
Given any $V_i , V_j \in V(R)$, let $d_{i,j}$ denote the density of $(V_i,V_j)_{G'}$. 
Then
\begin{align}\label{first}
(c_1 -2d)|R| \leq \frac{(c_1-2d)nm}{m^2} \leq \frac{\sum_{v\in V_i} \left( d^+_{G'}(v)-|V_0|\right)}{m^2}\le
\sum _{V_j \in V(R)} d_{i,j}
\end{align}
by Lemma~\ref{dilemma}. 
Thus 
\begin{align*}
d^+ _{R_m} (V_i) & = \sum _{V_j \in V(R_m)} \left\lfloor \frac{d_{i,j}}{\beta} \right\rfloor \geq \frac{1}{\beta} 
\sum _{V_j \in V(R)} d_{i,j} -|R_m| 
\stackrel{(\ref{first})}{\ge} (c_1-2d -\beta)\frac{|R_m|}{\beta}\\
& >(c_1-3d)\frac{|R_m|}{\beta}.
\end{align*}
So indeed $\delta ^+ (R_m) >(c_1-3d)|R_m|/\beta$. Similar arguments can be used to show
that $\delta ^- (R_m) >(c_1-3d)|R_m|/\beta$
and $\Delta ^0 (R_m)< (c_2 +2 \eps)|R_m|/\beta$.%
    \COMMENT{Indeed, given any $V_i \in V(R)$, $ c_2 nm\ge \sum_{v\in V_i}d^+_{G'}(v)
\ge m^2\sum _{V_j \in V(R)} d_{i,j}$. So $\sum _{V_j \in V(R)} d_{i,j}< (c_2 +2 \eps )|R|$. Thus 
$d^+ _{R_m} (V_i)= \sum _{V_j \in V(R_m)} \lfloor \frac{d_{i,j}}{\beta} \rfloor
\leq \frac{1}{\beta} \sum _{V_j \in V(R)} d_{i,j}< (c_2 +2 \eps)|R_m|/\beta$.}
\endproof

We will also need the well-known fact that for any cycle $C$ of the reduced multigraph $R_m$
we can delete a small number of vertices from the clusters in~$C$ in order to ensure that each edge
of $C$ corresponds to a super-regular pair. We include a proof for completeness.

\begin{lemma}\label{superreg}
Let $C=V_{j_1}\dots V_{j_s}$ be a cycle in the reduced multigraph~$R_m$ as in Lemma~\ref{multimin}. For each
$t=1,\dots,s$ let $(V_{j_t}V_{j_{t+1}})_{k_t}$
denote the edge of $C$ which joins $V_{j_t}$ to~$V_{j_{t+1}}$ (where $V_{j_{s+1}}:=V_{j_1}$). Then
we can choose subclusters $V'_{j_t}\subseteq V_{j_t}$ of size $m':=(1-4\eps)m$ such
that $(V'_{j_t},V'_{j_{t+1}})_{S_{j_t,j_{t+1},k_t}}$
is $(10\eps,\beta)$-super-regular (for each $t=1,\dots,s$).
\end{lemma}
\proof
Recall that for each $t=1,\dots, s$ the digraph $S_{j_t,j_{t+1},k_t}$ corresponding to the
edge $(V_{j_t}V_{j_{t+1}})_{k_t}$ of $C$ is $\eps$-regular and has density $\beta\pm \eps$.
So $V_{j_t}$ contains at most $2\eps m$ vertices whose outdegree in $S_{j_t,j_{t+1},k_t}$
is either at most $(\beta-2\eps)m$ or at least $(\beta+2\eps)m$. Similarly, there are
at most $2\eps m$ vertices in $V_{j_t}$ whose indegree in $S_{j_{t-1},j_{t},k_{t-1}}$
is either at most $(\beta-2\eps)m$ or at least $(\beta+2\eps)m$. Let $V'_{j_t}$ be a
set of size $m'$ obtained from $V_{j_t}$ by deleting all these vertices (and some additional
vertices if necessary). It is easy to check that $V'_{j_1},\dots,V'_{j_t}$ are subclusters
as required.
\endproof

Finally, we will use the following crude version of the  
fact that every $[\eps,\eps']$-regular pair contains a subgraph of 
given maximum degree $\Delta$ whose average degree is close to $\Delta$.
\begin{lemma} \label{boundmax}
Suppose that $0<1/n \ll \eps', \eps \ll d_0 \le d_1 \ll 1$ and that $(A,B)$ is an $[\eps,\eps']$-regular pair
of density $d_1$ with $n$ vertices in each class.
Then $(A,B)$ contains a subgraph $H$ whose maximum degree is at most $d_0n$
and whose average degree is at least $d_0n/8$.
\end{lemma}
\proof
Let $A'' \subseteq A$ be the set of vertices of degree at least $2d_1n$ and define $B''$ similarly.
Then $|A''|,|B''| \le \eps n$. Let $A':=A \setminus A''$ and $B':=B \setminus B''$.
Then $(A',B')$ is still $[2\eps,2\eps']$-regular of density at least $d_1/2$. 
Now consider a spanning subgraph $H$ of $(A',B')$ which is obtained from $(A',B')$
by including each edge with probability $d_0/3d_1$. So the expected degree
of every vertex is at most $2d_0n/3$ and the expected number of edges of 
$H$ is at least $d_0 (n-\eps n)^2/6$. Now apply the Chernoff bound on the binomial distribution
in Proposition~\ref{chernoff}
to each of the vertex degrees and to the total number of edges in $H$
to see that with high probability $H$ has the desired properties.\COMMENT{If a vertex in $(A',B')$ already has degree $\leq d_0 n$ then don't need
to apply Chernoff as they already have desired property. Apply Chernoff for the remaining vertices (i.e. doing this $\leq 2n \ll e^n$ times) and the 
edge set. Note that we are using $\mathbb P(d_H (x) \geq (1+\eps)2d_0n/3 ) \leq
 \mathbb P (|d_{H} (x)- \mathbb E(d_H (x))| \geq \eps \mathbb E (d_H (x)))
\leq 2 e^{-\frac{\eps^2}{3} d_0 ^2 n/(3d_1)} \leq 2e^ {-c n}$ for some constant $c>0$ (since $\mathbb E (d_H (x)) \geq d_0 ^2 n/(3d_1)$).
Further $\mathbb P (e (H) \leq (1-\eps)d_0 (n-\eps n )^2 /6) \leq \mathbb P (e (H) \leq (1-\eps) \mathbb E (e(H)))\leq 
P (|e(H)- \mathbb E (e(H))| \geq \eps \mathbb E (e(H)) \leq
2e^{-cn^2}$. (Since $\mathbb E (e(H)) \geq d_0 (n-\eps n)^2/6$.)}
\endproof

\section{Useful results} \label{4}
\subsection{$1$-factors in multidigraphs}

Our main aim in this subsection is to show that the reduced multidigraph
$R_m$ contains a collection of `almost' 1-factors
which together cover almost all the edges of~$R_m$ (see Lemma~\ref{multifactor1}). To prove this
we will need the following result which implies $R_m$ contains many edges between any two sufficiently
large sets. The second part of the lemma will be used in Section~\ref{sec:shifted}.

\begin{lemma}\label{keevashmult}
Let $0 <1/n  \ll 1/M'\ll \eps \ll \beta \ll \eta \ll d \ll c,d'\ll 1$. Suppose that $G$ is an oriented graph
of order $n$ with $\delta ^0 (G)\geq (1/2-\eta)n$. Let $R$ and $R_m$ denote the reduced digraph
and the reduced multidigraph of $G$ obtained by applying Lemma~\ref{dilemma}
(with parameters $\eps,d, M'$ and $\eps, \beta ,d, M'$ respectively). Let $L:=|R|=|R_m|$. 
Then the following properties hold.
\begin{itemize}
\item[{\rm (i)}] Let $X \subseteq V(R_m)$ be such that
$\delta ^0 (R_m [X]) \geq (1/2-c)|X|/\beta$.
Then for all (not necessarily disjoint) subsets $A$ and $B$ of $X$ of size at
least $(1/2-c)|X|$ there are at least $|X|^2/(60 \beta )$ directed edges from $A$ to $B$ in~$R_m$.
\item[{\rm (ii)}] Let $R'$ denote the spanning subdigraph of $R$ obtained by deleting all edges which correspond
to pairs of density at most $d'$ (in the pure digraph $G'$). Then $\delta^0(R')\ge (1/2-2d')L$
and for all (not necessarily disjoint) subsets $A$ and $B$ of $V(R')$ of size at
least $(1/2-c)L$ there are at least $L^2/60$ directed edges from $A$ to $B$ in~$R'$.
\end{itemize}
\end{lemma}
\proof We first prove~(i). 
Recall that for every edge $V_iV_j$ of $R$ there are precisely $\lfloor d_{i,j} /\beta \rfloor$
edges from $V_i $ to $V_j$ in $R_m$, where $d_{i,j}$ denotes the density of $(V_i, V_j )_{G'}$.
But $d_{i,j}+d_{j,i} \leq 1 $ since $G$ is oriented and so $R_m$ contains at most
$1/\beta$ edges between $V_i$ and $V_j$ (here we count the edges in both directions).

By deleting vertices from $A$ and $B$ if necessary we may assume that $|A|=|B|=(1/2-c)|X|$.
We will distinguish two cases. Suppose first that $|A\cap B|>|X|/5$ and let $Y:=A\cap B$. 
Define $\overline{Y}:=X \backslash Y$ and $\overline{A \cup B}:=X\backslash (A\cup B)$.
Then
\begin{eqnarray*}
2e(A,B)& \ge & 2e(Y)=\sum_{V\in Y} d_{R_m [X]}(V)-e(Y,\overline{Y})-e(\overline{Y},Y)\\
& {\ge} & |Y|(1-2c)|X|/\beta -|Y|(|X|-|Y|)/\beta
=|Y|(|Y|-2c|X|)/\beta\ge |X|^2/(30\beta).
\end{eqnarray*}
So suppose next that $|A\cap B|\le |X|/5$. Then $|\overline{A\cup B}|\le |X|-|A|-|B|+|A\cap B|\le (1/5+2c)|X|$.
Therefore,
\begin{eqnarray*}
e(A,B) & \ge & \sum_{V\in A} d^+_{R_m [X]}(V)-e(A,\overline{A\cup B})-e(A)\\
& {\ge} &
|A|(1/2-c)|X|/\beta-|A||\overline{A\cup B}|/\beta-|A|^2/(2\beta)\\
 & \ge & |A|[(1/2-c)-(1/5+2c)-(1/2-c)/2]|X|/\beta\ge |X|^2/(60\beta),
\end{eqnarray*}
as required. 

To prove~(ii) we consider the weighted digraph $R'_w$ obtained from $R'$ by giving each edge
$V_iV_j$ of $R'$ weight $d_{i,j}$. Given a cluster $V_i$, we write $w^+(V_i)$ for the sum
of the weights of all edges sent out by~$V_i$ in~$R'_w$. We define $w^-(V_i)$ similarly and
write $w^0(R'_w)$ for the minimum of $\min\{w^+(V_i),w^-(V_i)\}$ over all clusters~$V_i$.
Note that $\delta^0(R')\ge w^0(R'_w)$. Moreover, Lemma~\ref{dilemma} implies that
$d^{\pm} _{G'\backslash V_0} (x) > (1/2 -2d)n \text{ \  for all  } x \in V(G' \backslash V_0)$.
Thus each $V_i \in V(R')$ satisfies
$$
(1/2-2d)nm \leq e_{G'}(V_i, V(G')\backslash V_0) \leq m^2 w^+ (V_i) +(d'm^2)L
$$
and so $w^+(V_i) \geq (1/2 -2d-d')L>(1/2-2d')L$. Arguing in the same way for inweights gives us
$\delta^0(R')\ge w^0 (R'_w) > (1/2-2d')L$. Let $A,B\subseteq V(R')$ be as in~(ii).
Similarly as in~(i) (setting $\beta:=1$ and $X:=V(R')$ in the
calculations) one can show that the sum of all weights of the edges from~$A$ to $B$ in~$R'_w$
is at least $L^2/60$. But this implies that $R'$ contains at least $L^2/60$ edges from~$A$ to~$B$.
\endproof

\begin{lemma}\label{multifactor1}
Let $0 <1/n  \ll 1/M'\ll \eps \ll \beta \ll \eta \ll d \ll c \ll 1$. Suppose that $G$ is an oriented
graph of order $n$ with $\delta ^0 (G)\geq (1/2-\eta)n$. Let
$R_m$ denote the reduced multidigraph of $G$ with parameters $\eps, \beta ,d$ and $M'$
obtained by applying Lemma~\ref{dilemma}. Let $r:= (1/2-c) |R_m|/\beta$.
Then there exist edge-disjoint collections $\mathcal F_1, \dots , \mathcal F_r$ of vertex-disjoint
cycles in $R_m$ such that each $\mathcal F_i$ covers all but 
at most $c|R_m|$ of the clusters in $R_m$.
\end{lemma}
\proof
Let $L:=|R_m|$.
Since $\Delta ^0 (G) \leq n-\delta ^0 (G) \leq (1/2+\eta )n$, Lemma~\ref{multimin} implies that
\begin{align}\label{multideg} \delta ^0 (R_m) \geq (1/2-4d)\frac{L}{\beta} \text{ \ \ and \ \ }\Delta ^0 (R_m) \leq (1/2+2\eta)\frac{L}{\beta}.
\end{align}

First we find a set of clusters $X \subseteq V(R)$ with the following properties:
\begin{itemize}
\item $|X|= cL$,
\item $|N^{\pm} _{R_m} (V_i) \cap X| = (1/2\pm 5d)\frac{cL}{\beta}$ for all $V_i \in V(R_m)$. 
\end{itemize}
We obtain $X$ by choosing a set of $cL$ clusters uniformly at random.
Then each cluster $V_i$ satisfies
$$ \mathbb E (|N^{\pm} _{R_m} (V_i) \cap X|) = c |N^{\pm} _{R_m} (V_i)|
\stackrel{(\ref{multideg})}{=} c(1/2\pm 4 d)\frac{L}{\beta}.$$ 
Proposition~\ref{chernoff} for the hypergeometric distribution now implies that
with nonzero probability $X$ satisfies our desired conditions. (Recall that $N^{+} _{R_m} (V_i)$ is a multiset. Formally Proposition~\ref{chernoff}
does not apply to multisets. However, for each $j=1, \dots , 1/\beta$ we can apply Proposition~\ref{chernoff} to the set of all those clusters which
appear at least $j$ times in $N^+ _{R_m} (V_i)$, and similarly for $N^- _{R_m}(V_i)$.)\COMMENT{We have $2L$ sets $N^{\pm} _{R_m} (V_i)$. For each 
set we apply Chernoff at most $1/ \beta $ times. So we use Chernoff $\leq 2L/\beta \ll e^L$ times. Let $N^+ _j (V_i)$ denote the set of all clusters 
that appear $\geq j$ times in $N^+ _{R_m} (V_i)$ (similarly define $N^- _j (V_i)$). If  $|N^{\pm} _j (V_i)|\leq \eps cL$ don't need to apply Chernoff.
Otherwise set $X^+ _{ij} :=X\cap N^+ _j (V_i)$. Then 
$\mathbb P (X^+ _{ij} \geq (1+\eps)c |N^+ _j (V_i)| \text{  or  } X^+ _{ij} \leq (1-\eps)c |N^+ _j (V_i)|) \leq
\mathbb P(|X^+ _{ij} -\mathbb E (X^+ _{ij})| \geq \eps \mathbb E(X^+ _{ij}) )\leq 2e^{-\frac{\eps ^3}{3} c^2 L}$. So
whp $|N^+ _{R_m} (V_i)\cap X| =(1 \pm \eps )c |N^+ _{R_m} (V_i)| \pm \eps c L / \beta =(1/2\pm 5d)\frac{cL}{\beta}$ as desired.}

Note that
$$  d^{\pm} _{R_m \backslash X} (V_i) = \left( \frac{1}{2}- \frac{c}{2}\pm 5d\right)\frac{L}{\beta}$$
for each $V_i \in V(R_m\backslash X)$.
We now add a small number of \emph{temporary edges} to $R_m \backslash X$ in order to turn it
into an $r'$-regular multidigraph where 
$r':=( \frac{1}{2}- \frac{c}{2}+5d)\frac{L}{\beta}$.
We do this as follows. As long as $R_m \backslash X$ is not $r'$-regular there
exist $V_i , V_j \in V(R_m \backslash X)$  such that $V_i$ has outdegree less than $r'$ and
$V_j$ has indegree less than $r'$. In this case we add an edge from $V_i$ to $V_j$. (Note we may have
$i=j$, in which case  we add a loop.) 

We decompose the edge set of $R_m \backslash X $ into $r'$ 1-factors
$\mathcal F'_1, \dots, \mathcal F'_{r'}$.
(To see that we can do this, consider the bipartite multigraph $H$ where both vertex classes $A,B$ consist of a copy of
$V(R_m \backslash X)$ and we have $s$ edges between $a \in A$ and $b \in B$ if
there are precisely $s$ edges from $a$ to~$b$ in $R_m \backslash X$, 
including the temporary edges. Then $H$ is regular and so has a perfect matching.%
   \COMMENT{We need Hall for bipartite multigraphs here. But the proof of Hall still works.}
This corresponds to a $1$-factor $\mathcal F'_1$. Now remove the edges of $\mathcal F'_1$ from~$H$
and continue to find $\mathcal F'_2,\dots,\mathcal F'_{r'}$ in the same way.)
Since at each cluster we added at most
$20d \frac{L}{\beta}$ temporary edges, all but at most $20 \sqrt{d}\frac{L}{\beta}$ of the $\mathcal F'_i$
contain at most $\sqrt{d} L$ temporary edges. By relabeling if necessary we may assume that
$\mathcal F'_1, \dots , \mathcal F'_{r}$ are such $1$-factors.
We now remove the temporary edges from each of these $1$-factors, though we still refer to
the digraphs obtained in this way as $\mathcal F'_1, \dots , \mathcal F'_{r}$.
So each $\mathcal F'_i$ spans $R_m \backslash X$ and consists of cycles and at most $\sqrt {d} L$ paths.

Our aim is to use the clusters in $X$ to piece up these paths into cycles in order to obtain
edge-disjoint directed subgraphs $\mathcal F_1, \dots , \mathcal F_{r}$ of $R_m$ where each
$\mathcal F_i$ is a collection of vertex-disjoint cycles and
$\mathcal F'_i \subseteq \mathcal F_i$.%
  \COMMENT{note we could have 2-cycles.}

Let $P'_1, \dots, P'_\ell$ denote all the paths lying in one of
$\mathcal F'_1, \dots , \mathcal F'_{r}$ (so $\ell \leq  \sqrt{d}L r\leq  \sqrt {d} L^2/\beta$).
Our next task is to find edge-disjoint paths and cycles $P_1, \dots, P_\ell$ of 
length~$5$ in $R_m$ with the following properties.
\begin{itemize}
\item[(i)] If $P'_j$ consists of a single cluster $V_{j'} \in V(R)$ then $P_j$ is a cycle consisting of $4$ clusters in $X$ as well as $V_{j'}$.
\item[(ii)] If $P'_j$ is a path of length $\geq 1$ then $P_j$ is a path whose
startpoint is the endpoint of $P'_j$. Similarly the endpoint of $P_j$ is the startpoint of $P'_j$.
\item[(iii)] If $P'_j$ is a path of length $\ge 1$ then the internal clusters in the path $P_j$ lie in~$X$.
\item[(iv)] If $P'_{j_1}$ and $P' _{j_2}$ lie in the same $\mathcal F'_i$ then $P_{j_1}$ and $P _{j_2}$ are vertex-disjoint. 
\end{itemize}
So conditions (i)--(iii) imply that $P'_j \cup P_j$ is a directed cycle for each $1 \leq j \leq \ell$.
Assuming we have found such paths and cycles $P_1, \dots, P_\ell$, we define
$\mathcal F_1 , \dots , \mathcal F_{r}$ as follows. Suppose $P'_{j_1}, \dots ,P'_{j_t}$ are the paths in
$\mathcal F'_i$. Then we obtain $\mathcal F_i$ from $\mathcal F'_i$ by adding the paths and cycles 
$P_{j_1}, \dots ,P_{j_t}$  to $\mathcal F'_i$. Condition~(iv) ensures that the $\mathcal F_i$
are indeed collections of vertex-disjoint cycles.

It remains to show the existence of $P_1, \dots, P_\ell$. Suppose that for some $j\le \ell$
we have already found $P_1, \dots, P_{j-1}$ and now need to define $P_j$.
Consider $P'_j$ and suppose it lies in $\mathcal F'_i$. Let $V_{a}$ denote the startpoint of $P'_j$ and
$V_{b}$ its endpoint.

We call an edge $(V_{i_1}V_{i_2})_k$ in $R_m$ {\emph{free}} if it has not been used
in one of $P_1, \dots ,P_{j-1}$. Let $B$ be the set of all those clusters $V\in X$
for which at least $c|X|/\beta$ of the edges at $V$ in $R_m [X]$ are not free.
Our next aim is to show that $B$ is small. More precisely,  
$$|B|\leq d^{1/4} L.$$
To see this, note that $3(j-1)\leq 3\ell\leq 3 \sqrt{d}\frac{L^2}{\beta}$ edges of $R_m[X]$ lie in one of
$P_1, \dots , P_{j-1}$. Thus,
$2\cdot 3\sqrt{d}\frac{L^2}{\beta} \geq \frac{c|X|}{\beta}|B|=\frac{c^2L|B|}{\beta}$.
(The extra factor of~2 comes from the fact that we may have counted edges at the vertices in~$B$ twice.)
Since $c\gg d$ this implies that $|B|\leq d^{1/4} L$,
as desired. We will only use clusters in $X':=X\backslash B$ when constructing $P_j$. 
Note that $V_{a}$ receives at most $|B|/\beta \leq d^{1/4}L/\beta$ edges from $B$
in $R_m$.

Since we added at most $20dL/\beta $ temporary edges to $R_m \backslash X$ per cluster, $V_{a}$ can be the
startpoint or endpoint of at most $20dL/\beta $ of the paths $P'_1, \dots, P'_{j-1}$.
Thus $V_{a}$ lies in at most $20dL/\beta$ of the paths and cycles
$P_1, \dots , P_{j-1}$. In particular, at most $40dL/\beta$ edges at $V_a $ in $R_m$ are not free.%
    \COMMENT{40 instead of 20 since each $P_s$ containing $V_a$ could be a cycle}
We will avoid such edges when constructing $P_j$. 

For each of $P_1, \dots , P_{j-1}$ we have used $4$ clusters in $X$. 
Let $P'_{j_1}, \dots , P'_{j_t}$ denote the paths which lie in $\mathcal F'_i$ (so $t \leq \sqrt{d}L$).
Thus at most $4\sqrt{d}L$ clusters in $X$ already lie in the paths and cycles $P_{j_1}, \dots ,P_{j_t}$.
So for $P_j$ to satisfy (iv), the inneighbour of $V_{a}$ on $P_j$ must not be one of these clusters.
Note that $V_{a}$ receives at most 
$4\sqrt{d}L/\beta$ edges in $R_m$ from these clusters.

Thus in total we cannot use $d^{1/4}L/\beta+40dL/\beta +4\sqrt{d}L/\beta
\leq 2d^{1/4}L/\beta$ of the edges which $V_{a}$ receives from $X$ in $R_m$. But
$|N_ {R_m} ^- (V_{a}) \cap X|\ge (\frac{1}{2}-5d)cL/\beta\gg 2d^{1/4}L/\beta$ and so
we can still choose a suitable cluster $V_{a^-}$ in
$N_ {R_m} ^- (V_{a}) \cap X$ which will play the role of the inneighbour of $V_{a}$ on $P_j$. 
Let $(V_{a^-}V_{a})_{k_5}$ denote the corresponding free edge in $R_m$ which we will use in $P_j$.

A similar argument shows that we can find a cluster $V_{b^+} \not = V_{a^-}$ to play the role
of the outneighbour of $V_{b}$ on $P_j$. So $V_{b^+} \in X'$, $V_{b^+}$ does not lie on any of
$P_{j_1}, \dots, P_{j_t}$ and there is a free edge $(V_{b}V_{b^+})_{k_1}$ in $R_m$.

We need to choose the outneighbour $V_{b^{++}}$ of $V_{b^+}$ on $P_j$ such that
$V_{b^{++}}\in X'\setminus \{V_{a^-}\}$, $V_{b^{++}}$ has not been used in $P_{j_1},\dots,P_{j_t}$
and there is a free edge from $V_{b^+}$ to $V_{b^{++}}$ in $R_m$. Let $A_1$ denote the set of all 
clusters in $X'$ which satisfy these conditions. Since $V_{b^+}\in X'$ at most $c|X|/\beta $
edges at $V_{b^+}$ in $R_m [X]$ are not free. So $V_{b^+}$ sends out at least 
$(1/2-5d)\frac{|X|}{\beta}-c\frac{|X|}{\beta}-\frac{|B\cup \{V_{a^-}\}|}{\beta}\geq (1/2-2c)\frac{|X|}{\beta}$
free edges to $X'\setminus \{V_{a^-}\}$ in $R_m$. On the other hand, as before one can show that $V_{b^+}$ sends at most
$4\sqrt{d}L/\beta$ edges to clusters in $X'$ which already lie in $P_{j_1}, \dots, P_{j_t}$. Hence,
$|A_1|\geq \beta[(1/2-2c)|X|/\beta- 4\sqrt{d}L/\beta]\geq (1/2-3c)|X|$. 

Similarly we need to choose the inneighbour $V_{a^{--}}$ of $V_{a^-}$ on $P_j$
such that $V_{a^{--}}\in X'\setminus \{V_{b+}\}$, $V_{a^{--}}$ has not been used in 
$P_{j_1}, \dots, P_{j_t}$ and so that $R_m$ contains a free edge from $V_{a^{--}}$ to $V_{a^-}$.
Let $A_2$ denote the set of all clusters in $X'$ which satisfy these conditions.
As before one can show that $|A_2| \geq (1/2-3c)|X|$.

Recall that $\delta ^0 (R_m [X])\geq (1/2-5d)|X|/\beta$ by our choice of~$X$.
Thus Lemma~\ref{keevashmult}(i) implies that $R_m [X]$ contains at least $|X|^2/(60\beta)= c^2L^2/(60 \beta)$ edges from $A_1$
to $A_2$. Since all but at most $5\ell\le 5\sqrt{d}L^2/\beta$ edges of $R_m$ are free,
there is a free edge $(V_{b^{++}}V_{a^{--}})_{k_3}$ from $A_1$ to $A_2$.
Let $(V_{b^+}V_{b^{++}})_{k_2}$ be a free edge from $V_{b^+}$ to $V_{b^{++}}$ in $R_m$
and let $(V_{a^{--}}V_{a^{-}})_{k_4}$ be a free edge from $V_{a^{--}}$ to $V_{a^{-}}$
(such edges exist by definition of $A_1$ and $A_2$).
We take $P_j$ to be the directed path or cycle which consists of the edges 
$(V_{b}V_{b^+})_{k_1}$, $(V_{b^{+}}V_{b^{++}})_{k_2}$, $(V_{b^{++}} V_{a^{--}})_{k_3}$,
$(V_{a^{--}} V_{a^-})_{k_4}$ and $(V_{a^-}V_a)_{k_5}$.
\endproof


\subsection{Spanning subgraphs of super-regular pairs}

Frieze and Krivelevich~\cite{fk} showed that every $(\eps,\beta)$-super-regular pair $\Gamma$
contains a regular subgraph $\Gamma'$ whose density is almost the same as that of $\Gamma$.
The following lemma is an extension of this, where we can require $\Gamma'$ to have a given
degree sequence, as long as this degree sequence is almost regular.

\begin{lemma}\label{fandk}
Let $0 < 1/m\ll\eps  \ll \beta \ll \alpha ' \ll \alpha \ll 1$. Suppose that $\Gamma =(U,V)$ is an
$(\eps , \beta+\eps)$-super-regular pair where $|U|=|V|=m$. Define $\tau:= (1- \alpha)\beta m$. Suppose we have a
non-negative integer $x_i \leq \alpha ' \beta m$ associated with each $u_i \in U$ and a
non-negative integer $y_i \leq \alpha ' \beta m$ associated with each $v_i \in V$ such that $\sum _{u_i 
\in U} x_i= \sum _{v_i \in V} y_i$. Then $\Gamma$ contains a spanning subgraph $\Gamma '$ in which
$c_i:= \tau-x_i$ is the degree of $u_i \in U$ and $d_i:= \tau-y_i$ is the degree of $v_i \in V$.
\end{lemma}
\proof We first obtain a directed network $N$ from $\Gamma$ by adding a source $s$ and a sink $t$. 
We add an edge $su_i$ of capacity $c_i$ for each $u_i \in U$ and an edge $v_i t$ of capacity $d_i$ for each
$v_i \in V$. We give all the edges in $\Gamma$ capacity $1$ and direct them
from $U$ to~$V$. 

Our aim is to show that the capacity of any cut is at least $\sum _{u_i \in U} c_i
= \sum _{v_i \in V} d_i$. By the max-flow min-cut theorem this would imply that $N$ admits
a flow of value $\sum _{u_i \in U} c_i$, which by
construction of $N$ implies the existence of our desired subgraph~$\Gamma'$.

So consider any $(s,t)$-cut $(S,\bar{S})$ where $S=\{s\} \cup S_1 \cup S_2$ with $S_1 \subseteq U$
and $S_2 \subseteq V$. Let $\bar S_1:=U \backslash S_1$ and $\bar S_2:=V \backslash S_2.$
The capacity of this cut is 
$$\sum _{u_i \in \bar S_1} c_i + \sum _{v_i \in S_2} d_i + e(S_1,\bar S_2) $$
and so our aim is to show that
\begin{align}\label{aim}
e(S_1,\bar S_2) \geq \sum _{u_i \in S_1} c_i - \sum _{v_i \in S_2} d_i.
\end{align}
Now 
\begin{align}\label{aim1}\sum _{u_i \in S_1} c_i - \sum _{v_i \in S_2} d_i \leq
|S_1|(1- \alpha)\beta m -|S_2|(1-\alpha -\alpha ')\beta m
\end{align} and similarly
\begin{align}\label{aim2}
\sum _{u_i \in S_1} c_i - \sum _{v_i \in S_2} d_i =\sum _{v_i \in \bar S_2} d_i -\sum _{u_i \in \bar S_1} c_i
\leq |\bar S_2|(1- \alpha)\beta m -|\bar S_1|(1-\alpha -\alpha ')\beta m.
\end{align}
By (\ref{aim1}) we may assume that $|S_1| \geq (1-2 \alpha ')|S_2|$. (Since otherwise $\sum _{u_i \in S_1} c_i - \sum _{v_i \in S_2} d_i <0$ and thus
(\ref{aim}) is satisfied.) Similarly by (\ref{aim2}) we may assume that $|\bar S_2| \geq (1-2 \alpha ')|\bar S_1|$.
Let $\alpha ^*:= \alpha '/ \alpha$.
We now consider several cases.

\medskip

\noindent
{\bf Case 1.} $|S_1|, |\bar S_2| \geq \eps m$ and $|S_1| \geq (1+\alpha^*)|S_2|.$

\smallskip \noindent Since $\Gamma$ is $(\eps, \beta +\eps)$-super-regular we have that 
\begin{align*}
e (S_1, \bar S_2) &\geq \beta |S_1|(m -|S_2|) \geq \beta m (|S_1|-|S_2|) \\
& = \left( |S_1|(1- \alpha)\beta m -|S_2|(1-\alpha -\alpha ')\beta m \right)+ \alpha \beta m |S_1|
-(\alpha +\alpha ')\beta m|S_2| \\
& \geq  |S_1|(1- \alpha)\beta m -|S_2|(1-\alpha -\alpha ')\beta m.
\end{align*}
(The last inequality follows since $\alpha  |S_1| \geq (\alpha +\alpha ')|S_2|$.) Together with~(\ref{aim1})
this implies (\ref{aim}).

\medskip

\noindent {\bf Case 2.} $|S_1|, |\bar S_2| \geq \eps m$, $|S_1| <(1+\alpha^*)|S_2|$ and
$|S_2| \leq (1- \alpha ^*)m.$

\smallskip \noindent Again since $\Gamma$ is $(\eps, \beta +\eps)$-super-regular we have that 
\begin{align}\label{eqaim}
e (S_1, \bar S_2) \geq \beta |S_1|(m -|S_2|)=\beta |S_1||\bar S_2|.
\end{align} 
As before, to prove~(\ref{aim}) we will show that 
$$e (S_1, \bar S_2) \geq  |S_1|(1- \alpha)\beta m -|S_2|(1-\alpha -\alpha ')\beta m.$$ Thus 
by (\ref{eqaim}) it suffices to show that
$\alpha m| S_1|- |S_1||S_2| +(1-\alpha-\alpha ')m|S_2| \geq 0$.
We know that
$|S_2| (1- \alpha -\alpha ') \geq |S_1|(1- \alpha - \alpha ^*)$ since $(1+\alpha^*)|S_2|>|S_1|.$%
\COMMENT{This follows as $1+\alpha ^* \leq (1- \alpha -\alpha ')/(1- \alpha - \alpha ^*)$ as 
$(1+\alpha ^*)(1- \alpha - \alpha ^*)=1- \alpha- \alpha \alpha ^* - (\alpha^*)^2 = 1- \alpha -\alpha ' - (\alpha^*)^2 \leq 1- \alpha -\alpha '.$}
Hence,
$\alpha |S_1|-|S_1|(1- \alpha ^*)+|S_2|(1-\alpha-\alpha ') \geq 0$. So
$\alpha m| S_1|- |S_1||S_2| +(1-\alpha-\alpha ')m|S_2| \geq 0$ as $|S_2| \leq (1- \alpha ^*)m.$
So indeed (\ref{aim}) is satisfied.

\medskip

\noindent {\bf Case 3.} $|S_1|, |\bar S_2| \geq \eps m$, $|S_1| <(1+\alpha^*)|S_2|$ and $|S_2| > (1- \alpha ^*)m.$

\smallskip \noindent
By~(\ref{aim2}) in order to prove~(\ref{aim}) it suffices to show that 
$$e (S_1, \bar S_2) \geq  |\bar S_2|(1- \alpha)\beta m -|\bar S_1|(1-\alpha -\alpha ')\beta m.$$
Since~(\ref{eqaim}) also holds in this case, this means that it suffices to show that
$\alpha | \bar S_2|m- |\bar S_1||\bar S_2| +(1-\alpha-\alpha ')|\bar S_1| m\geq 0$.%
\COMMENT{since we want to show
$| S_1||\bar S_2| \geq (1- \alpha ) |\bar S_2 |m-(1-\alpha-\alpha ')|\bar S_1| m$}
Since $|S_1|\geq (1- 2\alpha ')|S_2|$ and  $|S_2| > (1- \alpha ^*)m$ we have that $|S_1|> (1-\alpha )m$. 
Thus $\alpha |\bar S_2|m \geq |\bar S_1 ||\bar S_2|$ and 
so indeed  (\ref{aim}) holds.

\medskip

\noindent {\bf Case 4.} $|S_1| < \eps m \text{ and } |\bar S_2| \geq \eps m.$

\smallskip \noindent
Since $|S_1| \geq (1-2 \alpha ')|S_2|$ we have that $|S_2| \leq 2 \eps m$. Hence,
$$e(S_1, \bar S_2) \geq \beta m|S_1|-|S_1||S_2| \geq (\beta-2 \eps)m|S_1| \geq (1- \alpha)\beta m |S_1|$$
and so by (\ref{aim1}) we see that (\ref{aim}) is satisfied, as desired.

\medskip

\noindent {\bf Case 5.}  $|S_1| \ge \eps m \text{ and } |\bar S_2|< \eps m$.

\smallskip \noindent Similarly as in Case~4 it follows that
$e(S_1, \bar S_2) \geq  (1- \alpha)\beta m 
|\bar S_2|$
and so by (\ref{aim2}) we see that (\ref{aim}) is satisfied, as desired.

\medskip

\noindent
Note that we have considered all possible cases since
we cannot have that $|S_1|,|\bar S_2| < \eps m$. Indeed, if $|S_1|,|\bar S_2| < \eps m$ then $|S_2| \geq (1- \eps)m$ 
and as $|S_1|\geq (1- 2 \alpha ')|S_2|$ this implies $|S_1|\geq (1- 2 \alpha ')(1- \eps)m$, a contradiction.
\endproof


\subsection{Special $1$-factors in graphs and digraphs}

It is easy to see that every regular oriented graph $G$ contains a $1$-factor.
The following result states that  if $G$ is also dense, then (i) we can guarantee a $1$-factor
with few cycles. Such $1$-factors have the advantage that we can transform them into 
a Hamilton cycle by adding/deleting a comparatively small number of edges. 
(ii) implies that even if $G$ contains a sparse `bad' subgraph $H$, then there
will be a $1$-factor which does not contain `too many' edges of $H$.

\begin{lemma}\label{1factororiented}
Let $0<  \theta_1 ,\theta _2 , \theta _3 <1/2$ and $\theta _1 /\theta _3 \ll \theta _2$.
Let $G$ be a $\rho$-regular oriented graph whose order $n$ is sufficiently large and where 
$\rho:=\theta _3 n$. Suppose $A_1, \dots , A_{5n}$ are sets of vertices in $G$ with
$a_i:=|A_i|\ge n^{1/2}$. Let $H$ be an oriented subgraph of
$G$ such that $d^{\pm} _H (x) \leq \theta _1 n$ for all $x \in A_i$ (for each $i$).
Then $G$ has a $1$-factor $F$ such that
\begin{itemize}
\item[(i)] $F$ contains at most $n/(\log n)^{1/5}$  cycles;
\item[(ii)] For each $i$, at most $\theta _2 a_i$ edges of $H\cap F$  are incident to $A_i$.
\end{itemize}
\end{lemma}

To prove this result we will use ideas similar to those used by Frieze and Krivelevich~\cite{fk}.
In particular, we will use the following bounds on the number of perfect matchings in a bipartite graph.

\begin{thm}\label{matchingbounds}
Suppose that $B$ is a bipartite graph whose vertex classes have size $n$ and $d_1,\dots , d_n$ are the
degrees of the vertices in one of these vertex
classes. Let $\mu (B)$ denote the number of perfect matchings in $B$. Then
\begin{align*}
\mu (B) \leq \prod _{k=1} ^n (d_k !)^{1/d_k}.
\end{align*}
Furthermore, if $B$ is $\rho$-regular then 
\begin{align*}
 \mu (B) \ge\left(\frac{\rho}{n}\right) ^n n!.
\end{align*}
\end{thm}
The upper bound in Theorem~\ref{matchingbounds} was proved by Br\'egman~\cite{3}.
The lower bound is a consequence of the Van der Waerden conjecture which was
proved independently by Egorychev~\cite{eg} and Falikman~\cite{fa}. 

We will deduce (i) from the following result in~\cite{randmatch}, which in turn is similar to 
Lemma~2 in~\cite{fk}.

\begin{lemma} \label{2factor}
For all $\theta \le 1$ there exists $n_0=n_0(\theta)$ such that the following holds.
Let $B$ be a $\theta n$-regular bipartite graph whose vertex classes $U$ and $W$
satisfy $|U|=|W|=:n\ge n_0$. Let $M_1$
be any perfect matching from $U$ to $W$ which is disjoint from $B$. Let $M_2$
be a perfect matching chosen uniformly at random from the set
of all perfect matchings in $B$. Let $F=M_1 \cup M_2$ be the 
resulting $2$-factor. Then the probability that $F$ contains more than
$n/(\log n)^{1/5}$ cycles is at most $e^{-n}$. 
\end{lemma}

\medskip

\noindent
{\bf Proof of Lemma~\ref{1factororiented}.}
 Consider the $\rho$-regular bipartite graph $B$
whose vertex classes $V_1,V_2$ are copies of $V(G)$ and where $x \in V_1$ is joined to $y \in V_2$
if $xy$ is a directed edge in $G$. Note that
every perfect matching in $B$ corresponds to a $1$-factor of $G$ and vice versa.
Let $\mu (B)$ denote the number of perfect matchings of $B$. 
Then
\begin{align}\label{pmatch1}
\mu (B) \geq \left( \frac{\rho}{n} \right) ^n n!\geq \left( \frac{\rho}{n} \right) ^n \left( \frac{n}{e} \right) ^n = \left( \frac{\rho}{e} 
\right) ^n
\end{align}
by Theorem~\ref{matchingbounds}. Here we have also used Stirling's formula which implies that for
sufficiently large $m$,
\begin{align}\label{stirling}
\left( \frac{m}{e} \right) ^m \leq m! \leq \left( \frac{m}{e} \right) ^{m+1}. 
\end{align}
We now count the number $\mu _i (G)$ of $1$-factors of $G$ which contain more than $\theta _2 a_i$
edges of $H$ which are incident to $A_i$. Note that
\begin{align}\label{mured1} 
\mu _{i} (G) \leq \binom{2a_i}{\theta _2 a_i}(\theta _1 n)^{\theta _2 a_i} (\rho !)^ {(n-\theta _2 a_i)/\rho}.
\end{align}
Indeed, the term $\binom{2a_i}{\theta _2 a_i}(\theta _1 n)^{\theta _2 a_i}$ in~(\ref{mured1}) gives an
upper bound for the number of ways we can choose $\theta _2 a_i$ 
edges from $H$ which are incident to $A_i$ such that no two of these edges have the same
startpoint and no two of these edges have the same 
endpoint. The term $(\rho !)^ {(n-\theta _2 a_i)/\rho}$ in~(\ref{mured1}) uses the upper bound
in Theorem~\ref{matchingbounds} to give a bound on the number of $1$-factors in $G$
containing $\theta _2 a_i $ fixed edges. Now 
\begin{align}\label{rhobounda}
(\rho !)^{(n- \theta _2 a_i)/ \rho} \stackrel{(\ref{stirling})}{\leq}
\left( \frac{\rho}{e} \right) ^{(1+1/\rho)(n- \theta _2 a_i)} 
 \leq \left( \frac{\rho}{e} \right)^{n- \theta _2 a_i +1/\theta _3}
\end{align}
since $\rho = \theta _3 n$ and
\begin{align}\label{rhobound2a}
\left( \frac{e}{\rho} \right)^{ \theta _2 a_i -1/\theta _3} \leq 
\left( \frac{2e}{\theta _3 n} \right)^{\theta _2 a_i}
\end{align}
since $a_i \geq n^{1/2}$.
Furthermore,
\begin{align}\label{binomial1}
\binom{2a_i}{\theta _2 a_i} \leq  \frac{(2a_i)^{\theta _2 a_i} }{(\theta _2 a_i)!} 
\stackrel{(\ref{stirling})}{\leq} \left( \frac{2 e}{ \theta _2} \right) ^{\theta _2 a_i}.
\end{align}
So by (\ref{mured1}) we have that
\begin{align*}
\mu _{i} (G)  &\stackrel{(\ref{rhobounda}),(\ref{binomial1})}{\leq}
\left(\frac{2e}{\theta _2}\right)^ {\theta _2 a_i} (\theta _1 n)^{\theta _2 a_i}
\left( \frac{\rho}{e} \right)^{n- \theta _2 a_i + 1/\theta _3}\\
& \ \ \stackrel{(\ref{rhobound2a})}{\leq} \left( \frac{2e}{\theta _2} \theta _1 n
\frac{2 e}{\theta _3 n} \right)^{\theta _2 a_i}  \left( \frac{\rho}{e}\right) ^{n} 
\stackrel{(\ref{pmatch1})}{\leq}\left( \frac{4e^2 \theta _1}{\theta _2 \theta _3} \right)^{\theta _2 a_i} \mu (B)
\ll \frac{\mu(B)}{5n}
\end{align*}
since $\theta _1 / \theta _3 \ll \theta _2$, $ a_i \geq n^{1/2}$ and $n$ is sufficiently large.

Now we apply Lemma~\ref{2factor} to $B$ where $M_1$ is the identity matching (i.e.~every
vertex in $V_1$ is matched to its copy in $V_2$). Then a cycle of length $2\ell$ in $M_1 \cup M_2$
corresponds to a cycle of length $\ell$ in $G$.
So, since $n$ is sufficiently large, the number of $1$-factors of $G$ containing more than
$n/(\log n)^{1/5}$ cycles is at most $e^{-n} \mu (B)$. So there exists a $1$-factor $F$ of $G$ 
which satisfies (i) and (ii).
\endproof


\subsection{Rotation-Extension lemma}
The following lemma will be a useful tool when transforming $1$-factors into Hamilton cycles.
Given such a $1$-factor $F$, we will obtain a path $P$ by cutting up and connecting several cycles in $F$
(as described in the proof sketch in Section~\ref{sketch}).
We will then apply the lemma to obtain a cycle $C$ containing precisely the vertices of $P$.%
\begin{lemma}\label{rotationlemma}
Let $0 < 1/m \ll \eps \ll \gamma <1$. Let $G$ be an oriented graph on $n \geq 2m$ vertices.
Suppose that $U$ and $V$ are disjoint subsets of $V(G)$ of size $m$ with the following property:
\begin{align}\label{label}
\text{If }S \subseteq U, \ T \subseteq V \text{ are such that }|S|,|T| \geq \eps m
\text{ then }e_G (S,T) \geq \gamma |S||T|/2.
\end{align}
Suppose that $P=u_1 \dots u_k$ is a directed path in $G$ where $u_1 \in V$ and $u_k \in U$.
Let $X$ denote the set of inneighbours $u_i$ of $u_1$ which lie on $P$ so that $u_i \in U$
and $u_{i+1} \in V$. Similarly let $Y$ denote the set of outneighbours $u_i$ of $u_k$ which lie on
$P$ so that $u_i \in V$ and $u_{i-1} \in U$. Suppose that $|X|,|Y| \geq \gamma m$.
Then there exists a cycle $C$ in $G$ containing precisely the vertices of $P$
such that $|E(C)\backslash E(P)|\leq 5$. Furthermore, $E(P)\backslash E(C)$ consists of edges from~$X$
to~$X^+$ and edges from~$Y^-$ to~$Y$. (Here $X^+$ is the set of successors of vertices in~$X$
on $P$ and $Y^-$ is the set of predecessors of vertices in~$Y$ on~$P$.)
\end{lemma}
\proof
Clearly we may assume that $u_k u_1 \not \in E(G)$.
Let $X_1$ denote the set of the first $\gamma  m/2$ vertices in $X$ along $P$ and $X_2$ the set of the last 
$\gamma m/2$ vertices in $X$ along $P$. We define $Y_1$ and $Y_2$ analogously. 
So $X_1, X_2 \subseteq U$ and $Y_1, Y_2 \subseteq V$.
We have two cases to consider.

\medskip

\noindent {\bf{Case 1.}} All the vertices in $X_1$ precede those in $Y_2$ along $P$. 

\smallskip \noindent
Partition $X_1 = X_{11}  \cup X_{12}$ where $X_{11}$ denotes the set of the first
$\gamma m/4$ vertices in $X_1$ along $P$. We partition $Y_2$ into $Y_{21}$
and $Y_{22}$ analogously. Let $X_{12}^+$ denote the set of successors on $P$ of the vertices in $X_{12}$  and
$Y^-_{21}$ the set of predecessors of the vertices in $Y_{21}$. So
$X^+_{12} \subseteq V$ and $Y^-_{21} \subseteq U$. Further define
\begin{itemize}
\item $X'_{11} := \{ u_i \ | \ u_{i-1} \in X_{11} \text{ and } \exists \ \text{edge from $u_{i-1}$ to $X^+_{12}$} \} $ and
\item $Y'_{22} := \{ u_i \ | \ u_{i+1} \in Y_{22} \text{ and } \exists \ \text{edge from } Y^-_{21}  \text{ to }u_{i+1} \} $. 
\end{itemize}
So $X'_{11} \subseteq V$ and $Y'_{22} \subseteq U$. 

From (\ref{label}) it follows that $|X'_{11}|\geq  \frac{(\gamma /2)(\gamma m/4)|X^+_{12}|}{|X^+_{12}|} \geq  \eps m$
and similarly $|Y'_{22}| \geq \eps m$. Since $X'_{11} \subseteq V$ and $Y'_{22} \subseteq U$, by (\ref{label})
$G$ contains an edge $u_{i'}u_{i}$ from $ Y'_{22}$ to $X'_{11}$. 
Since $u_{i} \in X'_{11}$, by definition of $X'_{11}$ it follows that $G$ contains an edge $u_{i-1}u_{j}$ 
for some $u_{j} \in X^+_{12}$.
Likewise, since $u_{i'} \in Y'_{22}$, there is an edge $u_{j'} u_{i'+1}$
for some $u_{j'}\in Y^-_{21}$. Furthermore, $u_{j-1}u_1$ and  $u_k u_{j'+1}$ are edges of $G$ 
by definition of $X^+_{12}$ and $Y^-_{21}$. It is easy to check that the cycle
$$C=u_1 \dots u_{i-1} u_{j} u_{j+1} \dots u_{j'} u_{i'+1}u_{i'+2}\dots u_k
u_{j'+1} u_{j'+2} \dots u_{i'} u_{i}u_{i+1} \dots u_{j-1} u_1$$
has the required properties (see Figure~1). For example, $E(P)\backslash E(C)$ consists
of the edges $u_{i-1}u_i$, $u_{j-1}u_j$, $u_{j'}u_{j'+1}$ and $u_{i'}u_{i'+1}$. The former
two edges go from~$X$ to~$X^+$ and the latter two from~$Y^-$ to~$Y$.
\begin{figure}[htb!] \label{fig:rotation}
\begin{center}\footnotesize
\psfrag{1}[][]{\normalsize $u_1$}
\psfrag{2}[][]{\normalsize $u_{i-1}$}
\psfrag{3}[][]{\normalsize $u_{i}$}
\psfrag{4}[][]{\normalsize $u_{j-1}$}
\psfrag{5}[][]{\normalsize \ \ \ \ \ \ \ \ $u_{j} \in X_{12}^+$}
\psfrag{6}[][]{\normalsize $u_{j'}$}
\psfrag{7}[][]{\normalsize $u_{j'+1}$}
\psfrag{8}[][]{\normalsize $u_{i'}$}
\psfrag{9}[][]{\normalsize $u_{i'+1}$}
\psfrag{10}[][]{\normalsize $u_k$}
\includegraphics[width=0.70\columnwidth]{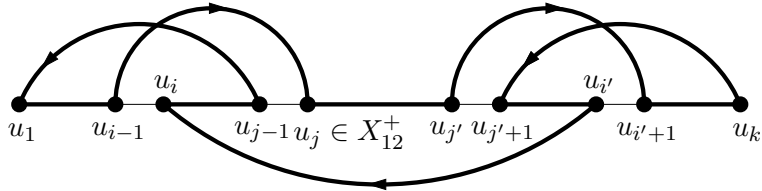}  
\caption{The cycle $C$ from Case~1}
\end{center} 
\end{figure}
 
\medskip

\noindent {\bf Case 2.} All the vertices in $Y_1$ precede those in $X_2$ along $P$.

\smallskip \noindent
Let $Y^- _1$ be the predecessors of the vertices in $Y_1$ and $X^+ _2$ the successors of the vertices in $X_2$ on $P$.
So $|Y^-_1|=|X^+_2|=\gamma m/2$ and $Y^- _1 \subseteq U$ and $X_2 ^+ \subseteq V$. Thus by~(\ref{label})
there exists an edge $u_{i}u_{j}\in E(G)$ from $Y^-_1$ to $X_2 ^+$. Again, it is easy to check
that the cycle
$$C=u_1 \dots u_{i} u_{j} u_{j+1} \dots u_k u_{i+1} u_{i+2} \dots u_{j-1} u_1$$
has the desired properties.
\endproof


\subsection{Shifted walks}\label{sec:shifted}
Suppose $R$ is a digraph and $F$ is a collection of vertex-disjoint cycles with $V(F) \subseteq V(R)$. 
A \emph{closed shifted walk $W$ in~$R$ with respect to $F$} is a walk in $R \cup F$ of the 
form 
$$W=c_1 ^+ C_1 c_1 c_2 ^+ C_2 c_2 \dots c^+ _{s-1} C_{s-1} c_{s-1} c^+ _s C_s c_s c_1 ^+,$$
where 
\begin{itemize}
\item $\{C_1, \dots , C_s\}$ is the set of all cycles in $F$; %
\item $c_i$ lies on $C_i$ and $c^+ _i$ is the successor of $c_i$ on $C_i$ for each $1 \le i \le s$;
\item $c_i c_{i+1} ^+$ is an edge of $R$ (here $c^+ _{s+1}:=c_1 ^+$).
\end{itemize}
Note that the cycles $C_1, \dots , C_s$ are not necessarily distinct. If a cycle $C_i$ in $F$ appears
exactly $t$ times in $W$ we say that $C_i$ is \emph{traversed $t$ times}. 
Note that a closed shifted walk $W$ has the property that for every cycle $C$ of $F$, every vertex of $C$ is
visited the same number of times by $W$.
The next lemma will be used in Section~\ref{merging} 
to combine cycles of $G$ which correspond to different cycles of $F$ into a single (Hamilton) cycle.
Shifted walks were introduced in~\cite{kelly}, where they were used for a similar purpose.

\begin{lemma}\label{shiftedwalk}
Let $0 <1/n  \ll 1/M'\ll \eps  \ll \eta \ll d \ll c\ll d' \ll 1$. Suppose that $G$ is an oriented
graph of order $n$ with $\delta ^0 (G)\geq (1/2-\eta)n$. Let
$R$ denote the reduced digraph of $G$ with parameters $\eps ,d$ and $M'$
obtained by applying Lemma~\ref{dilemma}. Let $L:=|R|$. Let $R'$ denote the spanning subgraph of $R$ obtained by 
deleting all edges which correspond to pairs of density at most $d'$ in the pure digraph~$G'$.
Let $F$ be  a collection of vertex-disjoint cycles with $V(F) \subseteq V(R')$ and $|V(F)| \ge (1-c)L$. 
Then $R'$ contains a closed shifted walk with respect to $F$ so that each cycle $C$ in $F$
is traversed at most $3L$ times.
\end{lemma}
\proof
Let $C_1, \dots , C_t$ denote the cycles of $F$.
We construct our closed shifted walk $W$ as follows: 
for each cycle $C_i$, choose an arbitrary vertex $a_i$ lying on $C_i$ and let $a_i^+$
denote its successor on~$C_i$. Let $U_i:=N^+_{R'}(a_i) \cap V(F)$ and
let $U_i^-$ be the set of predecessors of $U_i$ on $F$.
Similarly, let $V_i:=N^-_{R'}(a_i^+) \cap V(F)$ and let $V_i^+$ be the set of successors of $V_i$ on $F$. 
Since $\delta^0(R')\ge (1/2-2d')L$ by Lemma~\ref{keevashmult}(ii), we have
$|U_i^-|=|U_i| \ge (1/2-3d')L$ and $|V_i^+| =|V_i|\ge (1/2-3d')L$.
So by Lemma~\ref{keevashmult}(ii) there is an edge $u_i^-v_{i+1}^+$ from 
$U_i^-$ to $V_{i+1}^+$ in $R'$. Then we obtain a walk $W_i$ from $a^+_i$ to $a_{i+1}$ by first traversing
$C_i$  to reach $a_i$, then use the edge from $a_i$ to the successor $u_i$ of $u^-_i$, 
then traverse the cycle in~$F$ containing $u_i$ as far as $u_i^-$, then use the edge $u_i^-v_{i+1}^+$, then traverse
the cycle in~$F$ containing $v_{i+1}^+$ as far as $v_{i+1}$, and finally use the edge $v_{i+1}a_{i+1}$. (Here $a_{t+1}:=a_1$.)
$W$ is obtained by concatenating the $W_i$.
\endproof


\section{Proof of Theorem~\ref{main}} \label{5}
\subsection{Applying the Diregularity lemma}\label{applyDRL}
Without loss of generality we may assume that $0<\eta _1 \ll 1$. Define further constants satisfying 
\begin{align}\label{hier}0< 1/M'\ll \eps \ll \beta \ll \eta _2 \ll d\ll c \ll c' \ll
\gamma _1 \ll\gamma _2 \ll \gamma _3 \ll \gamma _4  \ll \gamma_5 \ll d' \ll \gamma  \ll \eta_1.
\end{align} 

Let $G$ be an oriented graph of order $n\gg M'$ such that $\delta ^0 (G) \geq (1/2-\eta _2)n$.
Apply the Diregularity lemma (Lemma~\ref{dilemma}) to $G$ with parameters $\eps, d$ and $M'$ to obtain
clusters $V_1, \dots,V_L$ of size $m$, an exceptional set $V_0$,
a pure digraph $G'$ and a reduced digraph $R$ (so $L=|R|$). 
Let $R'$ be the spanning subdigraph of~$R$
whose edges correspond to pairs of density at least~$d'$. So $V_iV_j$ is an edge of~$R'$ if
$(V_i,V_j)_{G'}$ has density at least~$d'$.

Let $R_m$ denote the reduced multidigraph of $G$ with parameters $\eps, \beta, d$ and $M'$. 
For each edge $V_iV_j$ of $R$ let $d_{i,j}$ denote the density of the $\eps$-regular pair $(V_i,V_j)_{G'}$.
Recall that each edge $(V_iV_j)_k \in E(R_m)$ is associated with the $k$th spanning subgraph
$S_{i,j,k}$ of $(V_i,V_j)_{G'}$ obtained by applying 
Lemma~\ref{split} with parameters $\eps, d_{i,j}$ and $K:=d_{i,j}/\beta$. 
Each $S_{i,j,k}$ is $\eps$-regular with density $\beta\pm \eps$.
Lemma~\ref{multimin} implies that
\begin{align}\label{Rmdeg}
\delta ^0 (R_m) \geq (1/2-4d)\frac{L}{\beta} \text{ \ \ and \ \ }\Delta ^0 (R_m) \leq (1/2+2\eta_2)\frac{L}{\beta}.
\end{align}
(The second inequality holds since $\Delta^0(G)\le n-\delta^0(G)\le (1/2+\eta_2)n$.)
Apply Lemma~\ref{multifactor1} to $R_m$ in order to obtain%
    \COMMENT{note we could choose more $1$-factors, but it is needed for 
later calculations that this $\gamma$ is `big'}
\begin{align}\label{rdef}
r:= (1-\gamma)L/2\beta
\end{align}
edge-disjoint collections $\mathcal F_1, \dots , \mathcal F_r$ of vertex-disjoint
cycles in $R_m$ such that each $\mathcal F_i$ contains all but 
at most $cL$ of the clusters in $R_m$.%
    \COMMENT{We could have $2$-cycles here. I don't think this will be a problem but if it is we can do a 
random split of clusters- double length of cycles...}
Let $V_{0,i}$ denote the set of all those vertices in $G$ which do not lie in clusters covered by~$\mathcal F_i$.
So $V_0 \subseteq V_{0,i}$ for all $1 \le i \le r$
and $|V_{0,i}|\leq |V_0|+cLm\leq (\eps+c)n$. 
We now apply Lemma~\ref{superreg} to each cycle in $\mathcal F_i$ to obtain subclusters of size
$m':=(1-4\eps)m$ such that the edges of $\mathcal F_i$ now correspond to $(10\eps,\beta)$-super-regular
pairs.  By removing one extra vertex from each cluster if necessary we may assume that $m'$ is even.
All vertices not belonging to the chosen subclusters of $\mathcal F_i$ are added to $V_{0,i}$. So now 
\begin{align}\label{v0}
|V_{0,i}|\leq 2cn.
\end{align}
We refer to the chosen subclusters as the clusters of $\mathcal F_i$ and still denote these
clusters by $V_1, \dots , V_L$. (This is a slight abuse of notation since the clusters of $\mathcal F_i$
might be different from those of $\mathcal F_{i'}$.)
Thus an edge $(V_{j_1}V_{j_2})_k$ in $\mathcal F_i$ corresponds to
the $(10\eps , \beta)$-super-regular pair $S' _{j_1, j_2 , k}:=(V_{j_1},V_{j_2})_{S_{j_1,j_2,k}}$.

Let $C_i$ denote the oriented subgraph of $G$ whose vertices are all those vertices belonging
to clusters in $\mathcal F_i$ such that for each $(V_{j_1}V_{j_2})_k \in E(\mathcal F_i)$ the edges between
$V_{j_1}$ and $V_{j_2}$ are precisely all the edges in $S' _{j_1,j_2,k}$.
Clearly $C_1, \dots , C_{r}$ are edge-disjoint.

We now define `random' edge-disjoint oriented subgraphs $H^+_1$, $H^-_1$, $H_2$, $H_{3,i}$, $H_4$ and
$H_{5,i}$ of $G$ (for each $i=1,\dots,r$).
$H^+_1$ and $H^-_1$ will be used in Section~\ref{sec:incorp} to incorporate the exceptional vertices
in $V_{0,i}$ into~$C_i$. $H_2$ will be used to choose the skeleton walks in Section~\ref{skel}.
The $H_{3,i}$ will be used in Section~\ref{4.6} to merge certain cycles. $H_4$ and the $H_{5,i}$ will
be used in Section~\ref{merging} to find our almost decomposition into Hamilton cycles.
We will choose these subgraphs to satisfy the following properties:

\medskip {\noindent\bf{Properties of $H^+_1$ and $H^-_1$.}}
\begin{itemize}
\item $H^+_1$ is a spanning oriented subgraph of $G$.
\item For all $x \in V(H^+_1)$, $\gamma _1 n \leq d^{\pm} _{H^+_1} (x) \leq 2 \gamma _1 n$.
\item For all $x \in V(H^+_1)$ and each $1 \leq i \leq r$, $|N ^{\pm} _{H^+_1} (x) \cap V_{0,i}|
\leq 4 \gamma _1 |V_{0,i}|$. 
\item $H^-_1$ satisfies analogous properties.
\end{itemize}

\medskip {\noindent\bf{Properties of $H_2$.}}%
\begin{itemize}
\item The vertex set of $H_2$ consists of precisely all those vertices of $G$ which lie in a cluster of $R$
(i.e.~$V(H_2)=V(G)\setminus V_0$). 
\item For each edge $(V_{j_1}V_{j_2})_k$ of $R_m$, $H_2$ contains a spanning oriented
subgraph of $S_{j_1, j_2, k}$ which forms an $\eps$-regular
pair of density at least $\gamma _2 \beta$.
\item All edges of $H_2$ belong to one of these $\eps$-regular pairs.
\item For all $x \in V(H_2)$, $d^{\pm}_{H_{2}} (x) \leq 2 \gamma_2 n$. 
\end{itemize}

\medskip {\noindent\bf{Properties of each $H_{3,i}$.}}
\begin{itemize}
\item The vertex set of $H_{3,i}$ consists of precisely all those vertices of $G$
which lie in a cluster of $\mathcal F_i$ (i.e.~$V(H_{3,i})=V(G)\setminus V_{0,i}$).
\item For each edge $(V_{j_1}V_{j_2})_k$ of $\mathcal F_i$, $H_{3,i}$ contains a spanning
oriented subgraph of $S' _{j_1, j_2, k}$ which forms a $(\sqrt{\eps}/2, 2\gamma _3 \beta )$-super-regular pair.
\item All edges in $H_{3,i}$ belong to one of these pairs.
\item Let $H_3$ denote the union of all the oriented graphs $H_{3,i}$. The last two properties
together with~(\ref{rdef}) imply that $d^{\pm}_{H_{3}}(x) \leq 3 \gamma_3 n$ for all $x \in V(H_3)$.
\end{itemize}

\medskip {\noindent\bf{Properties of $H_4$.}}
\begin{itemize}
\item The vertex set of $H_4$ consists of precisely all those vertices of $G$ which lie in a cluster of $R'$
(i.e.~$V(H_4)=V(G)\setminus V_0$). 
\item For each edge $V_{j_1}V_{j_2}$ of $R'$,  $(V_{j_1},V_{j_2})_{H_4}$ 
is $\eps$-regular of density at least $\gamma _4 d'$.
\item All edges in $H_4$ belong to one of these $\eps$-regular pairs.
\item For all $x \in V(H_4)$, $d^{\pm}_{H_{4}} (x) \leq 2 \gamma_4 n$. 
\end{itemize}

\medskip {\noindent\bf{Properties of each $H_{5,i}$.}}
\begin{itemize}
\item The vertex set of $H_{5,i}$ consists of precisely all those vertices of $G$ which lie in a
cluster of $\mathcal F_i$.
\item For each edge $(V_{j_1}V_{j_2})_k$ of $\mathcal F_i$, $H_{5,i}$ contains a spanning oriented
subgraph of $S' _{j_1, j_2, k}$ which forms a $(\sqrt{\eps}/2, 2\gamma _5 \beta )$-super-regular pair.
\item All edges in $H_{5,i}$ belong to one of these pairs.
\item Let $H_5$ denote the union of all the oriented graphs $H_{5,i}$. The last two properties together
with~(\ref{rdef}) imply that $d^{\pm}_{H_{5}}(x) \leq 3 \gamma_5 n$ for all $x \in V(H_5)$.
\end{itemize}

\medskip {\noindent\bf{Properties of each $S'_{i,j,k}$}.}
\begin{itemize}
\item For each edge $(V_{j_1}V_{j_2})_k$ of $\mathcal F_i$ the oriented subgraph obtained from
$S'_{j_1,j_2,k}$ by removing all the edges in $H^+_1,H^-_1,H_2,\dots,H_5$ is
$(\eps ^{1/3}, \beta _1)$-super-regular for some $\beta_1$ with%
   \COMMENT{Phrased in this way so that def of $\beta_1$ not too technical whilst still ensuring
can't have vertices of large degree in this super-reg pair-- I think this is important later on}
\begin{align}\label{beta1}
 (1-\gamma)\beta \leq \beta _1 \leq \beta.
\end{align}
\end{itemize}
The existence of $H^+_1$, $H^-_1$, $H_2$, $H_{3,i}$, $H_4$ and $H_{5,i}$ can be shown by considering suitable
random subgraphs of $G$ and applying the Chernoff bound in Proposition~\ref{chernoff}.
For example, to show that $H^+_1$ exists, consider a random subgraph of $G$ which
is obtained by including each edge of $G$ with probability $3\gamma_1$.%
   \COMMENT{Since the degrees in $G$ are about $n/2$ the expected degrees in~$H^+_1$ are about $3\gamma _1 n/2$}
Similarly, to define
$H_2$ choose every edge in $S_{j_1,j_2,k}$ with probability $3\gamma_2/2$ (for all $S_{j_1,j_2,k}$)
and argue as in the proof of Lemma~\ref{split}. 
Note that since $H_4$ only consists of edges between pairs of clusters $V_{j_1},V_{j_2}$ which
form an edge in $R'$, the oriented subgraphs obtained from the
$S'_{j_1,j_2,k}$ by deleting all the edges in $H^+_1,H^-_1,H_2,\dots,H_5$ may have densities which differ too much
from each other.
Indeed, if $V_{j_1}V_{j_2} \notin E(R')$, then the corresponding density will be larger.
However, for such pairs we can delete approximately a further $\gamma_4$-proportion of the edges 
to ensure this property
holds. Again, the deletion is done by considering a random subgraph obtained by deleting edges with probability $\gamma_4$.%

We now remove the edges in $H^+_1,H^-_1,H_2,\dots,H_5$ from each $C_i$.
We still refer to the subgraphs of $C_i$ and $S'_{j_1,j_2,k}$ thus obtained as $C_i$ and $S'_{j_1,j_2,k}$.


\subsection{Incorporating $V_{0,i}$ into $C_i$}\label{sec:incorp}
Our ultimate aim is to use each of the $C_i$ as a `framework' to piece together roughly
$\beta _1 m'$ Hamilton cycles in $G$. In this section we will incorporate the vertices in $V_{0,i}$,
together with some edges incident to these vertices, into $C_i$. For each $i=1,\dots,r$, let $G_i$
denote the oriented spanning subgraph of $G$ obtained from $C_i$ by adding the vertices of $V_{0,i}$.
So initially $G_i$ contains no edges with a start- or endpoint in $V_{0,i}$.
We now wish to add edges to $G_i$ so that
\begin{itemize}
\item[(i)] $d^{\pm} _{G_i} (x) \geq(1-\sqrt{c} ) \beta _1  m'$ where $x$ has neighbours only in $C_i$,
for all $x \in V_{0,i}$;
\item [(ii)] $|N^{\pm} _{G_i} (y) \cap V_{0,i}| \leq \sqrt{c} \beta _1 m'$ for all $y \in V(C_i)$;
\item[(iii)] $G_1, \dots ,G_{r}$ are edge-disjoint.
\end{itemize} 
For each $x \in V(G)$ we define $\mathcal L_x := \{ i \ | \ x \in V_{0,i}\}$ and let $L_x := |\mathcal L_x|$.
To satisfy~(i), we need to find roughly $L_x\beta_1 m'$ edges sent out by~$x$ (as well as $L_x\beta_1 m'$ edges received by~$x$)
such that none of these edges  already lies in any of the~$C_i$. It is not hard to check that such edges exist
(c.f.~(\ref{eq:freeedges}) below). However, if $L_x$ is small then there is not much choice to which $G_i$
with $i\in \mathcal{L}_x$ we add each of these edges and so it might not be possible to guarantee~(ii).
For this reason we reserved $H^+_1$ and $H^-_1$ in advance and for all
those $x$ for which $L_x$ is small we will use the edges at $x$ lying in these two graphs. More precisely, 
let $$B':= \left\{ x \in V(G) \ | \ L_ x \geq \frac{\gamma _1 n}{2\beta _1 m'}\right\}.$$ As indicated above, we now
consider the vertices in $B'$ and $V(G)\backslash B'$ separately.

First consider any $x\in V(G)\setminus B'$. Let $p:=2\beta _1 m'/ \gamma _1 n$ and
consider each edge $e$ sent out by $x$ in $H^+_1$. With probability $L_x p \leq 1$ we will assign $e$
to exactly one of the $G_i$ with $i \in \mathcal L_x$. More precisely, for each $i \in \mathcal L_x$
we assign $e$ to $G_i$ with probability $p$. So the probability $e$ is not assigned to any of the
$G_i$ is $1-L_x p \geq 0$. We randomly distribute the edges of~$H^-_1$ received by $x$ in an
analogous way amongst all the $G_i$ with $i\in \mathcal L_x$.

We proceed similarly for all
the vertices in $V(G)\setminus B'$, with the random choices being independent for different such
vertices. Since $H^+_1$ and $H^-_1$ are edge-disjoint from each other and from all the $C_i$,
the oriented graphs obtained from $G_1,\dots,G_r$ in this way will still be
edge-disjoint. Moreover, $\mathbb E (d^{\pm} _{G_i} (x))\geq \gamma _1 n p $
and $\mathbb E (d^{\pm} _{G_i [V_{0,i}]} (x)) \leq |V_{0,i}|p \leq 2 cn p$
for every $x \in V(G)\setminus B'$ and each $i\in \mathcal L_x$. Thus
\begin{align}\label{job}
\mathbb E(|N^{\pm}  _{G_i} (x) \cap V(C_i)|) \geq (\gamma _1  - 2c)n p \geq \beta _1 m'.
\end{align}
Let $B_i:= V_{0,i} \cap B'$ and $\bar{B}_i:=V_{0,i}\backslash B'$. Since
$|N^{\pm} _{H^+_1\cup H^-_1} (y) \cap V_{0,i}| \leq 8 \gamma _1 |V_{0,i}|$ for every $y \in V(C_i)$
(by definition of $H^+_1$ and $H^-_1$) we have that
\begin{align}\label{ctov} 
\mathbb E (|N^{\pm} _{G_i} (y) \cap \bar{B}_i| ) \leq 8\gamma _1 |V_{0,i}| p
\stackrel{(\ref{v0})}{\leq} 32 c \beta _1 m'.
\end{align}
Applying the Chernoff bound in Proposition~\ref{chernoff} (for the binomial distribution) for
each~$i$ and summing up the error probabilities for all~$i$ we see that with nonzero probability
the following properties hold:
\begin{itemize}
\item (\ref{job}) implies that $|N^{\pm} _{G_i} (x) \cap V(C_i)| \geq (1-\sqrt{c})\beta _1 m' $
for every $x \in \bar{B}_i$.
\item (\ref{ctov}) implies that $|N^{\pm} _{G_i} (y) \cap \bar{B}_i| \leq \sqrt{c}\beta_1 m'/2 $
for every $y \in V(C_i)$.
\end{itemize}
\COMMENT{Given a vertex $x \in V(G) \backslash B'$ and $i \in \mathcal L_x$ let $Z:= |N^+ _{G_i} (x) \cap V(C_i)|$.
So $\mathbb P (Z\geq (1- \eps )\beta _1 m' ) \leq \mathbb P (|Z-\mathbb E(Z)| \geq \eps \mathbb E(Z))\leq 2 e^{-\frac{\eps^2}{3 }\beta_1 m'}.$
Do this for every $x \in V(G) \backslash B'$ and $i \in \mathcal L_x$. So we use Chernoff at most $nr$ times for these cases. (We also repeat
for inneighbourhoods.) Note $4nr e^{-\frac{\eps ^2}{3} \beta _1 m'} \ll 1$, as desired.
Given any $i$ and $y \in V(C_i)$ if $|N^{-} _{H^+ _1} (y)\cap \bar{B}_i| \leq \sqrt{c} \beta _1 m'/2$ then clearly we don't need to apply Chernoff.
(Likewise if $|N^{-} _{H^+ _1} (y)\cap \bar{B}_i| \leq \sqrt{c} \beta _1 m'/2$.) Otherwise we apply Chernoff: Let $Z:=|N^- _{G_i} (y) \cap
\bar{B}_i|$. So $\mathbb P( Z\geq (1+\eps) 32c \beta _1 m') \leq \mathbb P (|Z-\mathbb E(Z)| \geq \eps \mathbb E (Z)) 
\leq 2 e^{-\frac{\eps^2}{6} \sqrt{c}\beta _1 m'}$. We apply Chernoff
for each such $y \in V(C_i)$ in each $G_i$. We also act analogously for $|N^+ _{G_i} (y) \cap
\bar{B}_i|$. Notice that $4 nr e^{-\frac{\eps^2}{6} \sqrt{c}\beta _1 m'}\ll 1$, so we have that whp our two conditions are satified.}
For each~$i$ we delete all the edges with both endpoints in~$V_{0,i}$ from~$G_i$.%
   \COMMENT{Have to do this already now since otherwise $(\ref{freedeg})$ and $(\ref{freedeg2})$
might not hold}

Having dealt with the vertices in $V(G)\setminus B'$, let us now consider any $x\in B'$.
We call each edge of $G$ with startpoint $x$ \emph{free} if it does not lie in any of
$C_i$, $H^+_1,H^-_1,H_2,\dots,H_5$ (for all $i=1,\dots,r$) and if the endpoint is not in $B'$.
Note that
$$ |B'| \frac{\gamma _1 n}{2 \beta_1 m'} \leq \sum ^{r} _{i=1} |V_{0,i}|
\stackrel{(\ref{v0})}{\leq} 2c r n \stackrel{(\ref{rdef})}{\leq} cn \frac{L}{\beta},$$ and so
$|B'| \leq \frac{2 c n}{\gamma _1}.$ 
So the number of free edges sent out by $x$ is at least
\begin{align}\label{eq:freeedges}
& ({1/2}-\eta _2)n -(\beta _1 + \eps^{1/3})m' (r-L_x)- 4\gamma _1 n -2\gamma _2 n -3\gamma _3 n -
2\gamma _4 n-3\gamma_5 n-|B'| \nonumber\\
& \stackrel{(\ref{rdef})}{\geq}({1/2}-\eta _2)n 
-(\beta +  \eps^{1/3})m'(1- \gamma )\frac{L}{2 \beta} + L_x \beta _1 m'-4\gamma _5 n -\frac{2 c n}{\gamma _1} \nonumber\\
& \stackrel{(\ref{hier})}{\geq} ({1/2}-\eta _2)n- \left( \frac{ \eps ^{1/3} n}{2 \beta}+\frac{n}{2} \right) +
\frac{\gamma  n}{4}+ L_x \beta _1 m' -5 \gamma _5 n
\stackrel{(\ref{hier})}{\geq} L_x \beta _1 m'.
\end{align}
We consider $L_x \beta _1 m'$ of these free edges sent out by $x$ and distribute them randomly amongst all
the $G_i$ with $i \in \mathcal L_x$. More precisely, each
such edge is assigned to $G_i$ with probability $1/L_x$ (for each $i \in \mathcal L_x$).
So for each $i \in \mathcal L_x$,
\begin{align}\label{freedeg}
\mathbb E (d^+ _{G_i} (x))= \beta_1 m'
\end{align}
and 
\begin{align}\label{freedeg2}
\mathbb E (d^+ _{G_i[V_{0,i}]} (x))\leq |V_{0,i}| \frac{1}{L_x} \stackrel{(\ref{v0})}{\leq}
2cn\left( \frac{2 \beta_1 m'}{\gamma_1 n}\right) =\frac{4c \beta_1 m'}{\gamma_1} \ll \sqrt{c} \beta _1 m'/4.
\end{align}
We can introduce an analogous definition of a free edge at $x$ but for edges whose endpoint is $x$.
As above we randomly distribute $L_x \beta _1 m'$ such edges amongst all the $G_i$ with $i \in \mathcal L_x$.
Thus for each $i \in \mathcal L_x$,
\begin{align}\label{freein}
\mathbb E (d^- _{G_i} (x))= \beta _1 m' \text{ \ \ and \ \ }
\mathbb E (d^- _{G_i[V_{0,i}]} (x)) \ll \sqrt{c} \beta _1 m'/4.
\end{align}
We proceed similarly for all vertices in $B'$, with the random choices being independent for
different vertices $x\in B'$. (Note that every edge of $G$ is free with respect to at most one vertex in~$B'$.)
Then using the lower bound on $L_x$ for all $x\in B'$ we have
\begin{align}\label{inC}
\mathbb E (|N^{\pm} _{G_i} (y) \cap B_i| ) \leq |V_{0,i}| \frac{2\beta_1 m'}{\gamma_1 n}
\stackrel{(\ref{v0})}{\leq} \sqrt{c} \beta _1 m' /4
\end{align}
for each $i=1,\dots,r$ and all $y \in V(C_i)$.
As before, applying the Chernoff type bound in Proposition~\ref{chernoff} for each $i$ and summing up
the failure probabilities over all~$i$ shows that with nonzero probability the following properties hold:
\begin{itemize}
\item (\ref{freedeg})--(\ref{freein}) imply that 
$|N^{\pm} _{G_i} (x) \cap V(C_i)|\geq (1- \sqrt{c}) \beta_1 m'  $ for each $x \in B_i$.
\item (\ref{inC}) implies that $|N^{\pm} _{G_i} (y) \cap B_i|\leq \sqrt{c} \beta_1 m' /2$
for each $y \in V(C_i)$.
\end{itemize}
\COMMENT{As before we get that whp the 1st condition is satisfied. Second condition is similar to before. Again we split into two cases:
If the number of free edges from vertices in $B_i$ to a vertex $y \in V(C_i)$ is $\leq \sqrt{c} \beta _1 m' /2$ then clearly we will have that
$|N^- _{G_i}(y) \cap B_i| \leq \sqrt{c} \beta _1 m' /2$. Otherwise we apply Chernoff as before.}
Together with the properties of $G_i$ established after choosing the edges at the vertices in $V(G)\setminus B'$
it follows that $|N^{\pm} _{G_i} (x) \cap V(C_i)|\geq (1- \sqrt{c}) \beta _1m' $ for every $x\in V_{0,i}$
and $|N^{\pm} _{G_i} (y) \cap V_{0,i}|\leq \sqrt{c} \beta _1 m' $ for every $y \in V(C_i)$. Furthermore,
$G_1,\dots,G_r$ are still edge-disjoint since when dealing with the vertices in $B'$ we only added
free edges. By discarding any edges assigned to $G_i$ which lie entirely in $V_{0,i}$
we can ensure that~(i) holds. So altogether (i)--(iii) are satisfied, as desired.


\subsection{Randomly splitting the $G_i$}\label{randomsplit}  
As mentioned in the previous section we will use each of the $G_i$ to piece together roughly
$\beta _1 m'$ Hamilton cycles of $G$. We will achieve this by firstly adding some more special
edges to each $G_i$ (see Section~\ref{skel}) and  then almost decomposing each $G_i$ into $1$-factors.
However, in order to use these $1$-factors to create Hamilton cycles we will need to ensure that
no $1$-factor contains a $2$-path with start- and endpoint in $V_{0,i}$, and midpoint in $C_i$.
Unfortunately $G_i$ might contain such paths. To avoid them, we will `randomly split' each~$G_i$.
 
We start by considering a random partition of each $V \in V(\mathcal F_i)$.
Using the Chernoff bound in Proposition~\ref{chernoff} for the hypergeometric distribution
one can show that there exists a partition of $V$ into subclusters $V'$ and $V''$ so that the
following conditions hold:
\begin{itemize}
\item $|V'|,|V''|=m'/2$ for each $V \in V(\mathcal F_i)$.
\item $|N^{\pm} _{G_i} (x) \cap \mathcal V'| \geq (1/2-\sqrt{c}) \beta _1 m'$ and
$|N^{\pm} _{G_i} (x) \cap \mathcal V''| \geq (1/2-\sqrt{c}) \beta _1 m'$ for each $x \in V_{0,i}$. (Here
$\mathcal V ':=\bigcup _{V \in V(\mathcal F_i)} V'$ and $\mathcal V '':=\bigcup _{V \in V(\mathcal F_i)} V''$.)
\end{itemize}
\COMMENT{Consider any $x \in V_{0,i}$ and some $V \in V (\mathcal F_i)$. Set $Z:=|N^+ _{G_i}(x) \cap V|$. If $Z < \eps m'/L$
then we won't apply Chernoff (note at most $\eps m'$ vertices in $G_i$ lie in such clusters $V$). So suppose that
$Z \geq \eps m'/L$. Let $Z':=|N^+ _{G_i} (x) \cap V'|$ and $Z''$ similarly. So $\mathbb E (Z')=\mathbb E (Z'')= Z/2 \geq \eps m'/(2L)$.
Thus $\mathbb P (Z' \leq (1-\eps) Z/2) \leq \mathbb P (|Z'-\mathbb E(Z')|\geq \eps \mathbb E(Z')) \leq 2 e^{- \frac{\eps^2}{3}\mathbb E(Z')}
\leq 2 e^{-\frac{\eps ^3}{6L}m'}$ (and similarly for $Z''$). We do this for each pair $x \in V_{0,i}$ and each $V \in V(\mathcal F_i)$.
Note that $|V_{0,i}| L \leq 2cnL$ and thus $8cnL e^{-\frac{\eps ^3}{6L}m'} \ll 1$. We deal with inneighbourhoods $N^- _{G_i} (x)$ analogously.
So whp we have that $|N^{\pm} _{G_i} (x) \cap \mathcal V'| \geq (1/2- \eps/2)[(1-\sqrt{c})\beta _1 m' - \eps m'] \geq (1/2- \sqrt{c} ) \beta _1 m'$
(and similarly for $\mathcal V''$).}
Recall that each edge $(V_{j_1}V_{j_2})_k \in E (\mathcal F_i)$ corresponds to
the $(\eps ^{1/3}, \beta _1)$-super-regular pair $S'_{j_1,j_2,k}$.
Let $\beta _2 := \beta _1/2$. So
\begin{align}\label{beta2}
(1/2-\gamma)\beta \stackrel{(\ref{beta1})}{\le} \beta _2 \stackrel{(\ref{beta1})}{\le} \beta /2.
\end{align}
Apply Lemma~\ref{split}(ii) to obtain a partition $E'_{j_1,j_2,k},E''_{j_1,j_2,k}$ of the edge set
of $S'_{j_1,j_2,k}$ so that the following condition holds:
\begin{itemize}
\item  The edges of $E'_{j_1,j_2,k}$ and $E''_{j_1,j_2,k}$ both induce an
$(\eps ^{1/4}, \beta _2)$-super-regular pair which spans $S'_{j_1,j_2,k}$.
\end{itemize}
We now partition $G_i$ into two oriented spanning subgraphs $G'_i$ and $G''_i$ as follows.
\begin{itemize}
\item The edge set of $G'_i$ is the union of all $E'_{j_1,j_2,k}$ (over all edges $(V_{j_1}V_{j_2})_k$
of $\mathcal F_i$) together with all the edges in $G_i$ from $V_{0,i}$ to $\mathcal V'$, and
all edges in~$G_i$ from $\mathcal V''$ to $V_{0,i}$.
\item The edge set of $G''_i$ is the union of all $E''_{j_1,j_2,k}$ (over all edges $(V_{j_1}V_{j_2})_k$
of $\mathcal F_i$) together with all the edges in $G_i$ from $V_{0,i}$ to $\mathcal V''$,
and all edges in~$G_i$ from $\mathcal V'$ to $V_{0,i}$. 
\end{itemize} 
Note that neither $G'_i$ nor $G''_i$ contains the type of $2$-paths we wish to avoid.
For each $i=1,\dots,r$ we use Lemma~\ref{split}(ii) to 
partition the edge set of each $H_{3,i}$ to obtain edge-disjoint
oriented spanning subgraphs $H'_{3,i}$ and $H''_{3,i}$ so that the following condition holds:
\begin{itemize}
\item For each edge $(V_{j_1}V_{j_2})_k$ in $\mathcal F_i$, both $H'_{3,i}$ and $H''_{3,i}$ contain a
spanning oriented subgraph of $S' _{j_1, j_2, k}$ which is $(\sqrt{\eps}, \gamma _3 \beta)$-super-regular. 
Moreover, all edges in $H'_{3,i}$ and $H''_{3,i}$ belong to one of these pairs.
\end{itemize} 
Similarly we partition the edge set of each $H_{5,i}$ to obtain edge-disjoint oriented
spanning subgraphs $H'_{5,i}$ and $H''_{5,i}$ so that the following condition holds:
\begin{itemize}
\item For each edge $(V_{j_1}V_{j_2})_k$ in $\mathcal F_i$, both $H'_{5,i}$ and $H''_{5,i}$ contain
a spanning oriented subgraph of $S' _{j_1, j_2, k}$ which is $(\sqrt{\eps}, \gamma _5 \beta)$-super-regular.
Moreover, all edges in $H'_{5,i}$ and $H''_{5,i}$ belong to one of these pairs.
\end{itemize} 
We pair $H'_{3,i}$ and $H'_{5,i}$ with $G'_i$ and pair $H''_{3,i}$ and $H''_{5,i}$ with $G''_i$.
We now have $2r$ edge-disjoint oriented subgraphs of $G$, namely $G'_1, G''_1, \dots , G'_r,G''_r$. 
To simplify notation, we relabel these oriented
graphs as $G_1, \dots, G_{r'}$ where
\begin{align}\label{r'}
r':=2r\stackrel{(\ref{rdef})}{=}(1-\gamma)L/\beta.
\end{align}
We similarly relabel the oriented graphs $H'_{3,1}, H''_{3,1}, \dots , H'_{3,r}, H''_{3,r}$ as 
$ H_{3,1}, \dots ,  H_{3,r'}$ and relabel $H'_{5,1}, H''_{5,1}, \dots , H'_{5,r}, H''_{5,r}$ as 
$ H_{5,1}, \dots ,  H_{5,r'}$ in such a way that $H_{3,i}$ and $H_{5,i}$ are the oriented graphs
which we paired with $G_i$. For each $i$ we still use the notation $\mathcal F_i$, $C_i$ and $V_{0,i}$ in the usual way.
Now (i) from Section~\ref{sec:incorp} becomes
\begin{itemize}
\item[(i$'$)] $d^{\pm} _{G_i} (x) \geq(1/2-\sqrt{c} ) \beta _1  m'$ where $x$ has neighbours only in $C_i$,
for all $x \in V_{0,i}$,
\end{itemize}
while (ii) and (iii) remain valid.

\subsection{Adding skeleton walks to the $G_i$}\label{skel}
Note that all vertices (including the vertices of $V_{0,i}$) 
in each $G_i$ now have in- and outdegree close to $\beta_2m'$.
In Section~\ref{nicefactor} our aim is to find a $\tau$-regular oriented subgraph of $G_i$, where 
\begin{align}\label{tau}
\tau:= (1-\gamma)\beta _2 m'.
\end{align} 
However, this may not be possible: suppose  for instance that $V_{0,i}$ consists of a 
single vertex $x$, $\mathcal F_i$ consists of 2 cycles $C$ and $C'$ and that
all outneighbours of $x$ lie on $C$ and all inneighbours lie on $C'$.
Then $G_i$ does not even contain a $1$-factor.
A similar problem arises if for example $V_{0,i}$ consists of a 
single vertex $x$, $\mathcal F_i$ consists of a single cycle $C=V_1 \dots V_t$,
all outneighbours of $x$ lie in the
cluster $V_2$ and all inneighbours in the cluster $V_8$.
Note that in both situations, the edges between $V_{0,i}$ and $C_i$ are not `well-distributed' or
`balanced'.
To overcome this problem, we add further edges to $C_i$ which will `balance out' the edges
between $C_i$ and $V_{0,i}$ which we added previously.
These edges will be part of the skeleton walks which we define below.
To motivate the definition of the skeleton walks it may be helpful to consider the second example 
above: Suppose that we add an edge $e$ from $V_1$ to $V_9$. Then $G_i$ now has a $1$-factor.
In general, we cannot find such an edge, but it will turn out that
we can find a collection of 5 edges fulfilling the same purpose.

A {\emph{skeleton walk}} $ S$ in $G$ with respect to $G_i$ is a collection of distinct
edges  $x_1x_2$, $x^-_2x_3$, $x_3^-x_4$, $x_4^-x_5$ and $x_5^-x_1$ of $G$  
with the following properties:
\begin{itemize}
\item $x_1 \in V_{0,i}$ and all vertices in $ V(S) \backslash \{ x_1 \}$ lie in $C_i$.
\item Given some $2\le j\le 5$, let $V\in V(\mathcal F_i)$ denote the cluster in~$\mathcal F_i$ containing
$x_j$ and let $C$ denote the cycle in $\mathcal F_i$ containing $V$.
Then $x_j^- \in V^-$, where $V^- $ is the predecessor of $V$ on $C$. 
\end{itemize}
Note that whenever $\mathcal S$ is a union of edge-disjoint skeleton walks and $V $ is a cluster
in $\mathcal F_i$, the number of edges in $\mathcal S$ whose endpoint is in $V$ is the same as
the number of edges in $\mathcal S$ whose startpoint is in $V^-$.
As indicated above, this `balanced' property will be crucial when finding a $\tau$-regular oriented subgraph of $G_i$ in
Section~\ref{nicefactor}. 

The 2nd, 3rd and 4th edge of each skeleton walk~$S$ with respect to~$G_i$ will lie in the
`random' graph~$H_2$ chosen in Section~\ref{applyDRL}. More precisely, each of these three edges
will lie in a `slice' $H_{2,i}$ of~$H_2$ assigned to~$G_i$. We will now partition~$H_2$
into these `slices' $H_{2,1},\dots,H_{2,r'}$. To do this,
recall that any edge $(V_{j_1}V_{j_2})_k$ in $R_m$ corresponds to an 
$\eps$-regular pair of density at least $\gamma _2 \beta$ in $H_2$. 
Here $V_{j_1}$ and $V_{j_2}$ are viewed as clusters in $R_m$, so $|V_{j_1}|=|V_{j_2}|=m$.
Apply Lemma~\ref{split}(i) to each such pair of clusters to find edge-disjoint oriented subgraphs
$H_{2,1}, \dots , H_{2,r'}$ of $H_2$ so that for each $H_{2,i}$ all the edges 
$(V_{j_1}V_{j_2})_k$ in $R_m$ correspond to $[\eps,5\beta \eps /L]$-regular pairs with density at
least $(\gamma _2 \beta -2\eps)\beta /L\geq \gamma _2 \beta ^2 /2L$ in $H_{2,i}$.

Recall that by~(i$'$) in Section~\ref{randomsplit} each vertex $ x \in V_{0,i}$ has at
least $ (1/2-\sqrt{c})\beta _1m' \geq \tau$
outneighbours in $C_i$ and at least $(1/2-\sqrt{c})\beta _1m'$ inneighbours in $C_i$. We pair
$\tau$ of these outneighbours $x^+$ with distinct inneighbours $x^-$. 
For each of these $\tau$ pairs $x^+  , x^-$ we wish to find a skeleton walk with respect to~$G_i$
whose $1$st edge is $xx^+$ and whose $5$th edge is  $x^-x$. We denote the union of these $\tau$ pairs $xx^+, x^-x$
of edges over all $x\in V_{0,i}$ by $\mathcal T_i$.

In Section~\ref{randomsplit} we partitioned each cluster $V \in V(\mathcal F_i)$ into subclusters $V'$ and~$V''$.
We next show how to choose the skeleton walks for all those $G_i$ for which each edge in $G_i$
with startpoint in $V_{0,i}$ has its endpoint in $\mathcal V'$ (and so each edge in $G_i$ with endpoint
in $V_{0,i}$ has startpoint in $\mathcal V''$).
The other case is similar, one only has to interchange $\mathcal V'$ and $\mathcal V''$.

\begin{claim}\label{skelG} 
We can find a set $\mathcal S_i$ of $\tau |V_{0,i}|$ skeleton walks
with respect to~$G_i$, one for each pair of edges in $\mathcal T_i$, such that $\mathcal S_i$
has the following properties:
\begin{itemize}
\item[(i)] For each skeleton walk in $\mathcal S_i$, its $2$nd, $3$rd and $4$th edge all
lie in $H_{2,i}$ and all these edges have their startpoint in $\mathcal V''$ and endpoint in $\mathcal V'$.
\item[(ii)] Any two of the skeleton walks in $\mathcal S_i$ are edge-disjoint.
\item[(iii)] Every $y \in V(C_i)$ is incident to at most $c^{1/5}\beta _2 m'$ edges
belonging to the skeleton walks in $\mathcal S_i$.
\end{itemize}
\end{claim}
Note that $|\mathcal S_i|=|\mathcal T_i|= \tau |V_{0,i}| \le 2c\beta _2 m'n$ by (\ref{v0}) and (\ref{tau}).
To find $\mathcal S_i$, we will first find so-called shadow skeleton walks
(here the internal edges are edges of $R_m$ instead of $G$).
More precisely, a \emph{shadow skeleton walk} $S'$ with respect to $G_i$ is a collection of
two edges $x_1x_2$, $x_5^-x_1$ of $G$ and three edges $(X_2^-X_3)_{k_2}$, $(X_3^-X_4)_{k_3}$,
$(X_4^-X_5)_{k_4}$ of $R_m$ with the following properties:
\begin{itemize}
\item  $x_1x_2$, $x_5^-x_1$ is a pair in $\mathcal T_i$.
\item $x_2 \in X_2$, $x_5^- \in X_5^-$ and each $X_j$ is a vertex of
a cycle in $\mathcal F_i$ and $X_j^-$ is the predecessor of $X_j$ on that cycle.
\end{itemize}
Note that in the second condition we slightly abused the notation: as $X_j$ is a cluster in~$R_m$,
it only corresponds to a cluster in~$\mathcal F_i$ (which has size~$m'$ and is a subcluster of
the one in~$R_m$). However, in order to simplify our exposition, we will use the same notation
for a cluster in~$R_m$ as for the cluster in~$\mathcal F_i$ corresponding to it.
 
We refer to the edge $(X^- _j X _{j+1})_{k_j}$ as the $j$th edge of the shadow skeleton walk $S'$.
Given a collection $\mathcal S'$ of shadow skeleton walks (with respect to $G_i$) we say an edge of $R_m$ is 
{\emph{bad}} if it is used at least $B:=c^{1/4} \beta^2 (m')^2/L$ times in $\mathcal S'$,
and {\emph{very bad}} if it is used at least $10B$ times in $\mathcal S'$.
We say an edge from $V$ to $U$ in $R_m$ is \emph{$(V,+)$-bad} if it is used at least $B$ times as a
$2$nd edge in the shadow skeleton walks of $\mathcal S'$. An edge from $W$ to $V$ in $R_m$ 
is \emph{$(V,-)$-bad} if it is used at least $B$ times as a $4$th edge in the shadow skeleton
walks of~$\mathcal S'$.

To prove Claim~\ref{skelG} we will first prove the following result.
\begin{claim}\label{shadow}
We can find a collection $\mathcal S'_i$ of $\tau |V_{0,i}|$ shadow skeleton walks with respect to $G_i$,
one for each of pair in $\mathcal T_i$, such that the following condition holds:
\begin{itemize}
\item For each $2 \le j \le 4$, every edge in $R_m$ is used at most $B$ times as a $j$th
edge of some shadow skeleton walk in $\mathcal S'_i$. In particular no edge in $R_m$ is very bad.
\end{itemize}
\end{claim} 
\proof
Suppose that we have already found $\ell < \tau |V_{0,i}|$ of our desired shadow skeleton walks for $G_i$. 
Let $xx^+, x^-x$ be a pair in $\mathcal T_i$ for which we have yet to define a
shadow skeleton walk. We will now 
find such a shadow skeleton walk $S'$. Suppose $x^+ \in V^+$ and $x^- \in W^-$, where
$V^+ , W^- \in V (\mathcal F_i )$. Let $V$ denote the predecessor of $V^+$ in $\mathcal 
F_i$ and $W$ the successor of $W^-$ in $\mathcal F_i$. 
We define $\mathcal V^+$ to consist of all those clusters $U \in V(\mathcal F_i)$ for which there
exists an edge from $V$ to $U$ in $R_m$ which is not $(V,+)$-bad.
By definition of $G_i$ (condition~(ii) in Section~\ref{sec:incorp}), each $y \in V(C_i)$
has at most $\sqrt{c} \beta _1 m'$ inneighbours in $V_{0,i}$ in $G_i$. So the number of $(V,+)$-bad edges is at most 
$$
\frac{\sqrt{c}\beta _1 (m')^2}{B}=\frac{\sqrt{c}\beta _1 (m')^2}{c^{1/4} \beta ^2   (m')^2/L} =
\frac{c^{1/4} \beta_1 L}{\beta^2} \stackrel{(\ref{beta1})}{\leq } \frac{c^{1/4}L}{\beta}.
$$ 
Together with~(\ref{Rmdeg}) this implies that
$$|\mathcal V^+|\geq (1/2 - 4d- c^{1/4})L \geq (1/2 -2c^{1/4})L.$$
Similarly we define $\mathcal W^-$ to consist of all those clusters $U \in V(\mathcal F_i)$ for which
there exists an edge from $U$ to $W$ in $R_m$ which is not $(W,-)$-bad.
Again, $|\mathcal W^-|\geq (1/2-2c^{1/4})L$. Let $\mathcal V$ denote the set
of those clusters which are the predecessors in $\mathcal F_i$ of a cluster in $\mathcal V^+$.
Similarly let $\mathcal W$ denote the set of those 
clusters which are the successors in $\mathcal F_i$ of a cluster in $\mathcal W^-$.
So $|\mathcal V|=|\mathcal V^+|$ and $|\mathcal W| = |\mathcal W^+|$. By Lemma~\ref{keevashmult}(i)
applied with $X=V(R_m)$ there exist at least $L^2/60 \beta$ edges in $R_m$ from $\mathcal V$ to $\mathcal W$. 
On the other hand, the
number of bad edges is at  most 
$$
\frac{3\tau |V_{0,i}|}{B}\stackrel{(\ref{v0}),(\ref{tau})}{\leq}
\frac{6\beta _2 m'cn}{c^{1/4} \beta ^2 (m')^2/L} 
\leq \frac{7 c^{3/4} \beta_2 L^2}{\beta^2 }
\stackrel{(\ref{beta2})}{\leq} \frac{7c^{3/4}L^2}{\beta}.
$$ 
So we can choose an edge $(XY)_k$ from $\mathcal V$ to $\mathcal W$ in $R_m$ which is
not bad. Let $X^+$ denote the successor of $X$ in $\mathcal F_i$  and $Y^-$ the predecessor of
$Y$ in $\mathcal F_i$. Thus $X^+ \in \mathcal V^+$ and $Y^- \in \mathcal W^-$ and so there is an
edge $(VX^+)_{k'}$ in $R_m$ which is not $(V,+)$-bad and an edge $(Y^-W)_{k''}$ which
is not $(W,-)$-bad. Let $S'$ be the shadow skeleton walk consisting of the edges $xx^+$, $(VX^+)_{k'}$,
$(XY)_k$, $(Y^-W)_{k''}$, and $x^-x$. Then we can add $S'$ to our collection of $\ell$ skeleton walks that
we have found already. 
\endproof 

We now use Claim~\ref{shadow} to prove Claim~\ref{skelG}.
{\removelastskip\penalty55\medskip\noindent{\bf Proof of Claim~\ref{skelG}.}} 
We apply Claim~\ref{shadow} to obtain a collection $\mathcal S'_i$ of shadow skeleton walks.
We will replace each edge of $R_m$ in these shadow skeleton walks with a distinct edge of $H_{2,i}$
to obtain our desired collection $\mathcal S_i$ of skeleton walks. 

Recall that each edge $(VW)_k$ in $R_m$ corresponds to an $[\eps, 5\eps\beta /L]$-regular pair
of density at least $\gamma _2 \beta ^2 /2L$ in $H_{2,i}$. 
Thus in $H_{2,i}$ the edges from $V''$ to $W'$ induce a $[3\eps, 10\eps\beta /L]$-regular pair%
   \COMMENT{3 instead of 2 since we first go to subclusters of size $m'$ and then partition each of these
into two}
of density $d_1 \geq \gamma _2 \beta ^2 /3L$. (Here $V',V''$ and $W',W''$ are the partitions of $V$ and
$W$ chosen in Section~\ref{randomsplit}.) 
Let $d_0:=80B/(m'/2)^2$ and note that $d_0 \le d_1$. So 
we can now apply Lemma~\ref{boundmax} to $(V'',W')_{H_{2,i}}$ to obtain a subgraph
$H'_{2,i}[V'',W']$ with maximum degree at most $d_0 m'/2 $
and at least $d_0 (m'/2)^2/8 = 10B $ edges.
We do this for all those edges in $R_m$ which are used in a shadow skeleton
walk in $\mathcal S'_i$. 

Since no edge in $R_m$ is very bad, for each $S' \in \mathcal S'_i$ we can replace an
edge $(VW)_k$ in $S'$ with a distinct edge $e$ from $V''$ to $W'$ lying in $H'_{2,i}[V'',W']$.
Thus we obtain a collection $\mathcal S_i$ of skeleton walks 
which satisfy properties (i) and (ii) of Claim~\ref{skelG}. Note that by the construction
of $\mathcal S_i$ every vertex $y \in V(C_i)$ is incident to at most
$d_0 m'L/(2\beta ) \ll c^{1/5} \beta _2 m'/2$
edges which play the role of a $2$nd, $3$rd or $4$th edge in a skeleton walk
in $\mathcal S_i$. Condition~(ii) in Section~\ref{sec:incorp} implies that $y$ is
incident to at most $2 \sqrt{c}\beta _1 m'$ edges which play the role of a 1st or 5th edge in a skeleton walk
in $\mathcal S_i$. So in total $y$ is incident to at most 
$c^{1/5} \beta _2 m'/2 +2 \sqrt{c}\beta _1 m' \leq c^{1/5} \beta _2 m'$ edges of the skeleton walks
in $\mathcal S_i$. Hence (iii) and thus the entire claim is satisfied. 
\endproof

We now add the edges of the skeleton walks in $\mathcal S_i$ to $G_i$. Moreover, for each $x\in V_{0,i}$
we delete all those edges at~$x$ which do not lie in a skeleton walk in~$\mathcal S_i$.


\subsection{Almost decomposing the $G_i$ into $1$-factors}\label{nicefactor}
Our aim in this section is to find a suitable collection of 1-factors in each~$G_i$ which together cover
almost all the edges of~$G_i$. In order to do this, we first choose a $\tau$-regular spanning oriented
subgraph $G^*_i$ of $G_i$ and then apply Lemma~\ref{1factororiented} to~$G^*_i$.

We will refer to all those edges in $G_i$ which lie in a skeleton walk in $\mathcal S_i$ as
\emph{red}, and all other edges in $G_i$ as {\emph{white}}. Given $V \in V(\mathcal F_i)$
and $x \in V$, we denote by $N^+ _{w} (x) $ the set of all those vertices which receive a white
edge from~$x$ in~$G_i$. Similarly we denote by $N^- _{w} (x) $ the set of all those vertices which send out a white
edge to $x$ in $G_i$. So $N^{+} _w (x)\subseteq V^+$ and $N^- _w (x) \subseteq V^-$, where $V^+$ and
$V^-$ are the successor and the predecessor of $V$ in $\mathcal F_i$.
Note that $G_i$ has the following properties:
\begin{itemize}
\item[$(\alpha _1)$] $d^{\pm} _{G_i} (x) = \tau $ for each $x \in V_{0,i}$. Moreover,
$x$ does not have any in- or outneighbours in $V_{0,i}$.
\item[$(\alpha _2)$] Every path in~$G_i$ consisting of two red edges has its midpoint in~$V_{0,i}$. 
\item[$(\alpha _3)$] For each $(V_jV^+ _j)_k\in E(\mathcal F_i)$ the white edges in $G_i$ from $V_j$
to $V^+ _j$ induce a $( \eps ^{1/4}, \beta _2)$-super-regular pair $(V_j , V^+ _j)_{G_i}$.
\item[$(\alpha _4)$] Every vertex $u \in V(C_i)$ receives at most $c^{1/5} \beta _2 m'$ red edges
and sends out at most $c^{1/5} \beta _2 m' $ red edges in $G_i$.
\item[$(\alpha _5)$] In total, the vertices in $G_i$ lying in a cluster $V_j \in V(\mathcal F_i)$
send out the same number of red edges as the vertices in $V^+ _j$ receive.
\end{itemize}
In order to find our $\tau$-regular spanning oriented subgraph of $G_i$, 
consider any edge $(V_jV^+ _j)_k \in E(\mathcal F_i)$. Given any 
$u_\ell \in V_j$, let $x_\ell$ denote the number of red edges sent out by $u_\ell$ in $G_i$.
Similarly given any $v_\ell \in V^+ _j$, let $y_\ell$ denote the number of red edges received by $v_\ell$
in $G_i$. By $(\alpha _4)$ we have that $x_\ell,y_\ell \leq c^{1/5} \beta _2 m'$ and 
by $(\alpha _5)$ we have that
$$\sum _{u_\ell \in V_j} x_\ell= \sum _{v_\ell \in V^+ _j} y_\ell.$$
Thus we can apply Lemma~\ref{fandk} to obtain an oriented spanning subgraph of
$(V_j , V^+ _j)_{G_i}$ in which each $u_\ell$ has outdegree $\tau-x_\ell$ and each $v_\ell$ has indegree
$\tau-y_\ell$. We apply Lemma~\ref{fandk} to each $(V_jV^+ _j)_k \in E(\mathcal F_i)$. The union of
all these oriented subgraphs together with the red edges in $G_i$ clearly yield a $\tau$-regular
oriented subgraph $G^* _i$ of $G_i$, as desired.

We will use the following claim to almost decompose $G^* _i$ into $1$-factors with certain useful properties.
\begin{claim}\label{nice1factors}
Let $G^*$ be a spanning $\rho$-regular oriented subgraph of $G_i$ where $\rho \geq \gamma \beta _2 m'$. 
Then $G^*$ contains a $1$-factor $F^*$ with the following properties:
\begin{itemize}
\item[{\rm(i)}] $F^*$ contains at most $n/(\log n)^{1/5}$ cycles.
\item[{\rm(ii)}] For each $V_j \in V(\mathcal F_i)$, $F^*$ contains at most $c'  m'$
red edges incident to vertices in $V_j$.
\item[{\rm(iii)}] Let $F^* _{red}$ denote the set of vertices which are incident to a red edge in $F^*$.
Then $|F^* _{red} \cap N^{\pm} _{H_{3,i}} (x)| \leq 2c'\gamma _3 \beta m'$ for each $x \in V(C_i)$.
\item[{\rm(iv)}] $|F^* _{red} \cap N^{\pm} _{w} (x)| \leq 2c' \beta _2 m'$ for each $x \in V(C_i)$.
\end{itemize}
\end{claim}
\proof A direct application of Lemma~\ref{1factororiented} to $G^*$ proves the claim. Indeed, we apply the lemma 
with $\theta _1 = (c^{1/5} \beta _2 m')/n $, $\theta _2 = c'$, $\theta _3 =\rho/n \geq (\gamma \beta _2 m')/n$ and
with the oriented spanning subgraph of $G^*$ whose edge set consists precisely of the red edges in $G^*$ playing
the role of $H$. Furthermore, the clusters in $V(\mathcal F_i)$ 
together with the sets $N^{\pm} _w (x)$ and $N^{\pm} _{H_{3,i}} (x)$ (for each $x \in V(C_i)$)
play the role of the $A_j$.
\endproof

Repeatedly applying Claim~\ref{nice1factors} we obtain edge-disjoint
$1$-factors $F_{i,1}, \dots , F_{i, \psi}$ of $G_i$
satisfying conditions (i)--(iv) of the claim, where 
\begin{align}\label{psi}
\psi := (1-2\gamma )\beta_2 m'.
\end{align}
Our aim is now to transform each of the $F_{i,j}$ into a Hamilton cycle using the edges of 
$H_{3,i}$, $H_4$
and $H_{5,i}$.


\subsection{Merging the cycles in $F_{i,j}$ into a bounded number of cycles}\label{4.6}
Let $D_1, \dots , D_{\xi}$ denote the cycles in $\mathcal F_i$ and define 
$V_G(D_k)$ to be the set of vertices in $G_i$ which lie in clusters in the cycle $D_k$. 
In this subsection, for each $i$ and $j$ we will merge the cycles in $F_{i,j}$ to obtain 
a $1$-factor $F'_{i,j}$ consisting of at most $\xi$ cycles.

Recall from Section~\ref{nicefactor} that we call the edges of $G_i$ which lie on a skeleton walk
in $\mathcal S_i$ red and the non-red edges of $G_i$ white. We call the edges of the `random' oriented
graph $H_{3,i}$ defined in Section~\ref{applyDRL} \emph{green}. (Recall that $H_{3,i}$ was modified in Section~\ref{randomsplit}.)
We will use the edges from $H_{3,i}$ to obtain $1$-factors
$F'_{i,1}, \dots , F' _{i , \psi}$ for each $G_i$ with the following properties:%
    \COMMENT{Previously $(\beta_4)$ was: Let $V \in V(\mathcal F_i)$ and $x \in V$.
If $x$ lies on a white edge $xy$ in
$F_{i,j}$ then $x$ either lies on the same edge $xy$ in $F'_{i,j}$ or $x$ lies on a green edge
$xz \in E(H_{3,i})$ in $F'_{i,j}$. In the latter case $y$ is
the endpoint of a green edge lying on the same cycle as $x$ in $F'_{i,j}$.}
\begin{itemize} 
\item[$(\beta _1)$] If $i\not = i'$ or $j \not = j'$ then $F'_{i,j}$ and $F'_{i',j'}$ are edge-disjoint.
\item[$(\beta _2)$] For each $V\in V( \mathcal F_i) $ all $x\in V$ which send out a white edge in
$F_{i,j}$ lie on the same cycle $C$ in $F'_{i,j}$.
\item[$(\beta _3)$] $|E(F'_{i,j})\backslash E(F_{i,j})|\le 6 n/(\log n)^{1/5}$ for all $i$ and $j$.
Moreover, $E(F'_{i,j})\backslash E(F_{i,j})$ consists of green and white edges only.
\item[$(\beta _4)$] For every edge in $F_{i,j}$ both endvertices lie on the same cycle in~$F'_{i,j}$.
\item[$(\beta _5)$] All the red edges in $F_{i,j}$ still lie in $F'_{i,j}$.
\end{itemize}
Before showing the existence of $1$-factors satisfying ($\beta_1$)--($\beta_5$), we will derive two
further properties ($\beta_6$) and ($\beta_7$) from them
which we will use in the next subsection.
So suppose that $F'_{i,j}$ is a $1$-factor satisfying the above conditions.
Consider any cluster $V\in V(\mathcal F_i)$. Claim~\ref{nice1factors}(ii)
implies that $F_{i,j}$ contains at most $c ' m'$ red edges with startpoint in $V$. So the cycle $C$
in $F'_{i,j}$ which contains all vertices $x\in V $ sending out a white edge in $F_{i,j}$
must contain at least $(1- c ')m'$ such vertices $x$. In particular there are at least $(1- c ')m'>c 'm'$
vertices $y \in V^+$ which lie on $C$. So some of these vertices $y$
send out a white edge in $F_{i,j}$. But by $(\beta _2)$ this means that $C$ contains all
those vertices $y \in V^+$ which send out a white edge in $F_{i,j}$. Repeating this argument shows
that $C$ contains all vertices in $V(D_k)$ which send out a white edge in $F_{i,j}$
(here $D_k$ is the cycle on $\mathcal F_i$ that contains $V$). Furthermore, by property $(\beta _4)$,
$C$ contains all vertices in $V(D_k)$ which receive a white edge in $F_{i,j}$. By
property~$(\alpha_2)$ in Section~\ref{nicefactor} no vertex of $C_i$ is both the a startpoint of a red edge in $G_i$
and an endpoint of a red edge in $G_i$. So this implies that all vertices in $V_G(D_k)$
lie on $C$. Thus if we obtain $1$-factors $F'_{i,1}, \dots , F'_{i,\psi}$ satisfying $(\beta _1)$--$(\beta _5)$
then the following conditions also hold:
\begin{itemize}
\item[$(\beta _6)$] For each $j=1,\dots,\psi$ and each $k=1,\dots,\xi$ all the vertices in $V_G(D_k)$ lie
on the same cycle in $F'_{i,j}$.
\item[$(\beta_7)$] For each $V\in V(\mathcal F_i)$ and each $j=1,\dots, \psi$ at most $c'm'$ vertices
in~$V$ lie on a red edge in~$F'_{i,j}$.  
\end{itemize}
(Condition~$(\beta_7)$ follows from Claim~\ref{nice1factors}(ii) and the `moreover' part of~$(\beta_3)$.)

For every $i$, we will define the 1-factors $F'_{i,1}, \dots , F' _{i , \psi}$
sequentially. Initially, we let $F'_{i,j}=F_{i,j}$. So the $F'_{i,j}$ satisfy
all conditions except $(\beta_2)$.
Next, we describe how to modify $F'_{i,1}$ so that it also satisfies $(\beta_2$).

Recall from Section~\ref{randomsplit} that for each edge $(VV^+)_k$ of $\mathcal F_i$ the pair
$(V,V^+)_{H_{3,i}}$ is $(\sqrt{\eps},\gamma_3\beta)$-super-regular and thus
$\delta^{\pm} (H_{3,i})  \geq (\gamma _3 \beta-\sqrt{\eps})  m' \ge \gamma _3 \beta m'/2$.
Furthermore,
whenever $V \in V(\mathcal F_i)$ and $x \in V$, the outneighbourhood of $x$ in $H_{3,i}$ lies in $V^+$
and the inneighbourhood of $x$ in $H_{3,i}$ lies in $V^-$.
Let $H'_{3,i}$ denote the oriented spanning subgraph of $H_{3,i}$ whose edge set consists of those
edges $xy$ of $H_{3,i}$ for which $x$ is not a startpoint of a red edge in our current $1$-factor $F'_{i,1}$
and $y$ is not an endpoint of a red edge in $F'_{i,1}$. 
Consider a white edge $xy$ in $F'_{i,1}$. Claim~\ref{nice1factors}(iii) implies that
$x$ sends out at most $2c' \gamma _3 \beta m'$ green edges $xz$ in $H_{3,i}$ which do not lie in $H'_{3,i}$.
So  $d^+ _{H'_{3,i}}(x)   \geq (1/2-2c') \gamma _3 \beta  m'$.
Similarly, $d^- _{H'_{3,i}}(y)   \geq (1/2-2c') \gamma _3 \beta  m'$.
(However, if $uv$ is a red edge in $F'_{i,1}$ then $d^+ _{H'_{3,i}}(u) =d^- _{H'_{3,i}}(v) =0$.)
Thus we have the following properties of $H_{3,i}$ and $H'_{3,i}$: 
\begin{itemize}
\item[($\gamma _1$)] For each $V\in V(\mathcal F_i)$ all the edges in~$H_{3,i}$ sent out by vertices
in~$V$ go to~$V^+$. 
\item[($\gamma _2$)] If $xy$ is a white edge in $F'_{i,1}$ then 
$ d^+ _{H'_{3,i}}(x) ,d^- _{H'_{3,i}}(y)\geq \gamma _3 \beta m'/3$.
\item[($\gamma _3$)] Consider any $V \in V(\mathcal F_i)$. Let $S \subseteq V$ and $T \subseteq V^+$
be such that $|S|,|T|\geq \sqrt{\eps} m' $. Then
$e_{H_{3,i}} (S,T) \geq \gamma _3 \beta |S||T|/2$.
\end{itemize}
If $F'_{i,1}$ does not satisfy ($\beta_2$), then it contains 
cycles $C\neq C^*$ such that  there is a cluster $V \in V(\mathcal F_i)$
and white edges $xy$ on~$C$ and $x^*y^*$ on~$C^*$
with $x,x^* \in V$ and $y,y ^* \in V^+$. 

We have 3 cases to consider. Firstly, we may have a green edge $xz \in E(H'_{3,i})$ such that $z$
lies on a cycle $C'\neq C$ in $F'_{i,1}$. Then $z \in V^+$ and $z$ is the endpoint
of a white edge in $F'_{i,1}$ (by $(\gamma_1)$ and the definition of $H'_{3,i}$).
Secondly, there may be a green edge $wy^* \in E(H'_{3,i})$ such that
$w$ lies on a cycle $C'\neq C^*$ in $F'_{i,1}$. So here
$w \in V$ is the startpoint of a white edge in $F'_{i,1}$. If neither of these cases hold,
then $N^+ _{H'_{3,i}} (x)$ lies on $C$ and $N^- _{H'_{3,i}}(y^*)$ lies on $C^*$. Since
$d^+ _{H'_{3,i}} (x) ,d^- _{H'_{3,i}}(y^*)\geq \gamma _3 \beta m'/3$ by~($\gamma _2$), we can use~($\gamma _3$)
to find a green edge $x'y'$ from $N^- _{H'_{3,i}}(y^*)$ to $N^+ _{H'_{3,i}} (x)$. Then $x' \in V$, $y' \in V^+$,
$x'$ is the startpoint of a white edge in $F'_{i,1}$ and $y'$ is the endpoint of a white edge
in $F'_{i,1}$. 

We will only consider the first of these 3 cases. The other cases can be dealt with analogously:
In the second case $w$ plays the role of $x$ and $y^*$ plays the role of $z$. In the third case $x'$
plays the role of $x$ and $y'$ plays the role of $z$.

So let us assume that the first case holds, i.e.~there is 
a green edge $xz \in E(H'_{3,i})$ such that $z$ lies on a cycle $C'\neq C$
in $F'_{i,1}$ and $z$ lies on a white edge $wz$ on $C'$. Let $P$ denote the directed path
$(C \cup C' \cup \{xz \} )\backslash \{ xy,wz \}$ from $y\in V^+$ to $w\in V$. 
Suppose that the endpoint $w$ of $P$ lies on a green edge $wv\in E(H'_{3,i})$ such that $v$ lies outside~$P$.
Then $v \in V^+$ is the endpoint of a white edge $uv$ lying on the cycle $C''$ in $F'_{i,1}$ which contains~$v$.
We extend $P$ by replacing $P$ and $C''$ with $(P \cup C'' \cup \{wv\}) \backslash \{ uv\}$. We make 
similar extensions if the startpoint $y$ of $P$ has an inneighbour in $H'_{3,i}$ outside $P$.
We repeat this `extension' procedure as long as we can. Let $P$ denote the path obtained in this way, say
$P$ joins $a \in V^+$ to $b \in V$. Note that $a$ must be the endpoint of a white edge in $F'_{i,1}$
and $b$ the startpoint of a white edge in $F'_{i,1}$.

We will now apply a `rotation' procedure to close $P$ into a cycle.
By ($\gamma _2$) $a$ has at least $\gamma _3 \beta m' /3$ inneighbours in $H'_{3,i}$ 
and $b$ has at least $\gamma _3 \beta m' /3$ outneighbours in $H'_{3,i}$ and all
these in- and outneighbours lie on $P$ since we could not extend $P$ any further. 
Let $X:=N^- _{H'_{3,i}}(a)$ and $Y:=N^+ _{H'_{3,i}}(b)$. So $|X|,|Y| \geq \gamma_3 \beta m' /3$ and
$X \subseteq V$ and $Y \subseteq V^+$ by~$(\gamma_1)$. Moreover, whenever $c\in X$ and $c^+$ is the
successor of~$c$ on~$P$, then either $cc^+$ was a white edge in $F'_{i,1}$ or $cc^+\in E(H'_{3,i})$.
Thus in both cases $c^+\in V^+$. So the set $X^+$ of successors in~$P$ of all the vertices in~$X$ lies in~$V^+$ 
and no vertex in~$X$ sends out a red edge in~$P$. Similarly one can show that the set $Y^-$ of predecessors in~$P$
of all the vertices in~$Y$ lies in~$V$ and no vertex in~$Y$ receives a red edge in~$P$.
Together with~($\gamma _3$) this shows that we can apply Lemma~\ref{rotationlemma} with
$P\cup H_{3,i}$ playing the role of $G$ and $V^+$ playing the role of~$V$ and $V$ playing the role of $U$
to obtain a cycle $\hat{C}$ containing precisely the vertices of $P$ such that $|E(\hat{C})\backslash E(P)|\leq 5$,
$E(\hat{C})\backslash E(P)\subseteq E(H_{3,i})$ and such that $E(P)\backslash E(\hat{C})$ consists of edges
from $X$ to~$X^+$ and edges from~$Y^-$ to~$Y$. Thus $E(P)\backslash E(\hat{C})$ contains no
red edges. Replacing $P$ with $\hat{C}$ gives us a $1$-factor 
(which we still call $F'_{i,j}$) with fewer cycles.
Also note that if the number of cycles is reduced by $\ell$, then 
we use at most $\ell+5 \le 6\ell$ edges in $H_{3,i}$ to achieve this.
So $F'_{i,j}$ still satisfies all requirements with the possible exception of ($\beta_2$).
If it still does not satisfy ($\beta_2$),
we will repeatedly apply this `rotation-extension' procedure until the current 
$1$-factor $F'_{i,1}$ also satisfies ($\beta_2$).
However, we need to be careful since we do not want to use edges of $H_{3,i}$ several times in this process.
Simply deleting the edges we use may not work as ($\gamma_2$) might fail later on 
(when we will repeat the above process for $F'_{i,j}$ with $j>1$).

So each time we modify $F'_{i,1}$, we also modify $H_{3,i}$ as follows. All the edges
from $H_{3,i}$ which are used in $F'_{i,1}$ are removed from $H_{3,i}$. All the edges which
are removed from $F'_{i,1}$ 
in the rotation-extension procedure are added to $H_{3,i}$. (Note that by~$(\beta_5)$ we never add red edges
to $H_{3,i}$.)  When we refer to
$H_{3,i}$, we always mean the `current' version of $H_{3,i}$, not the original one.
Furthermore, at every step we still refer to an edge of $H_{3,i}$ as green, even if initially the edge did
not lie in $H_{3,i}$. Similarly at every step we refer to the non-red edges of our current $1$-factor
as white, even if initially they belonged to $H_{3,i}$.

Note that if we added a green edge $xz$ into $F'_{i,1}$, then $x$
lost an outneighbour in $H_{3,i}$, namely $z$. However, $(\beta _5)$ implies
that we also moved some (white) edge $xy$ of $F'_{i,1}$ to $H_{3,i}$, where $y$ lies
in the same cluster $V^+ \in V(\mathcal F_i)$ as $z$ (here $x\in V$).
So we still have that $\delta^{+} (H_{3,i})  \geq \gamma _3 \beta  m'/3$. Similarly, at any stage
$\delta^{-} (H_{3,i})  \geq \gamma _3 \beta  m'/3$.
When $H_{3,i}$ is modified, then $H'_{3,i}$ is modified accordingly. This will occur
if we add some
white edges to $H_{3,i}$ whose start or endpoint lies on a red edge in $F'_{i,1}$. However,
Claim~\ref{nice1factors}(iv) implies that at any stage we still have 
$$ 
d^+_{H'_{3,i}}(x) ,d^-_{H'_{3,i}}(y)\geq (1/2-2c')\gamma _3 \beta m'- 2c' \beta _2 m'\geq\gamma _3\beta m'/3.
$$
Also note that by ($\beta_3$), the modified version of $H_{3,i}$ still satisfies
\begin{equation} \label{h3i}
e_{H_{3,i}} (S,T) \geq (\gamma _3 \beta-\sqrt{\eps}) |S||T| -6n/(\log n)^{1/5}
\geq \gamma _3 \beta |S||T|/2.
\end{equation}
So $H_{3,i}$ and $H'_{3,i}$ will satisfy ($\gamma_1$)--($\gamma_3$) throughout and thus the above argument
still works. So after at most $n/(\log n)^{1/5}$ steps $F'_{i,1}$ will also satisfy ($\beta_2$).

Suppose that for some $1< j \le \psi$ we have found $1$-factors $F'_{i,1}, \dots , F'_{i, j-1}$ 
satisfying $(\beta _1)$--$(\beta _5)$. 
We can now carry out the rotation-extension procedure for $F'_{i,j}$ in the same way as for 
$F'_{i,1}$ until $F'_{i,j}$ also satisfies ($\beta_2$). In the construction of $F'_{i,j}$, we
do not use the original $H_{3,i}$, but the modified version obtained in the
construction of $F'_{i,j-1}$. We then introduce the oriented spanning subgraph 
$H'_{3,i}$ of $H_{3,i}$ similarly as before (but with respect to the current $1$-factor
$F'_{i,j}$). Then all the above bounds on these graphs still hold, except that
in the middle expression of~(\ref{h3i}) we multiply the term $6n/(\log n)^{1/5}$ by $j$ to account
for the total number of edges removed from $H_{3,i}$ so far.
But this does not affect the next inequality.
So eventually, all the $F'_{i,j}$ will satisfy ($\beta_1$)--($\beta_5$).


\subsection{Merging the cycles in $F'_{i,j}$ to obtain Hamilton cycles}\label{merging}
Our final aim is to piece together the cycles in $F'_{i,j}$, for each $i$ and $j$, to
obtain edge-disjoint Hamilton cycles of $G$. Since
we have $\psi$ $1$-factors $F'_{i,1}, \dots , F'_{i,\psi}$ for each $G_i$, in total we will find
\begin{eqnarray*}
\psi r' & \stackrel{(\ref{r'}),(\ref{psi})}{=} & (1-2\gamma)\beta _2 m'(1-\gamma)L/\beta 
\stackrel{(\ref{beta2})}{\geq}(1-2\gamma)(1-\gamma)(1/2-\gamma)m' L\\
& \stackrel{(\ref{hier})}{\geq } &(1/2- \eta _1 )n
\end{eqnarray*}
edge-disjoint Hamilton cycles of $G$, as desired.

Recall that $R'$ was defined in Section~\ref{applyDRL}.
Given any $i$, apply Lemma~\ref{shiftedwalk} to obtain a closed shifted walk
$$W_i= U^+ _1 D'_1 U_1 U_2 ^+ D'_2 U_2 \dots U^+ _{s-1} D'_{s-1}U_{s-1} U^+ _s D' _s U_s U^+ _1$$ in $R'$
with respect to~$\mathcal F_i$  such that each cycle in $\mathcal F_i$ is traversed at most $3L$ times.
So $\{D'_1,\dots,D'_s\}$ is the set of all cycles in $\mathcal F_i$, $U^+ _k$ is the successor of $U_k$ on $D'_{k}$
and $U_k U_{k+1} ^+\in E(R')$ for each $k=1,\dots, s$ (where $U_{s+1}:=U_1$). Moreover,
\begin{align}\label{s}
s \leq 3L^2.
\end{align}
For each $1$-factor $F'_{i,j}$ we will now use the edges of $H_4$ and $H_{5,i}$ to obtain a Hamilton cycle $C_{i,j}$ with the following
properties:
\begin{itemize}
\item[(i)] If $i \not = i'$ or $j \not = j'$ then $C_{i,j}$ and $C_{i',j'}$ are edge-disjoint.
\item[(ii)] $E(C_{i,j})$ consists of edges from $F'_{i,j}$, $H_4$ and $H_{5,i}$ only.
\item[(iii)] There are at most $3L^2$ edges from $H_4$  lying in $C_{i,j}$.
\item[(iv)] There are at most $3L^2+5$ edges from $H_{5,i}$ lying in $C_{i,j}$.
\end{itemize}
For each $j$, we will use $W_{i}$ to `guide' us how to merge the cycles in $F'_{i,j}$ into the
Hamilton cycle $C_{i,j}$. Suppose that we have already defined $\ell < \psi r'$ of the Hamilton
cycles $C_{i',j'}$ satisfying (i)--(iv), but have yet to define $C_{i,j}$.
We remove all those edges which have been used in these $\ell$ Hamilton cycles from both~$H_4$
and $H_{5,i}$. 

For each $V \in V(\mathcal F_i)$, we denote by $V_w$ the subcluster of $V$ containing all those vertices
which do not lie on a red edge in $F'_{i,j}$. We refer to $V_w$ as the \emph{white subcluster of $V$}.
Thus $|V_w|\geq (1-c')m'$ by property~$(\beta_7)$ in Section~\ref{4.6}. Note that the outneighbours of
the vertices in $V_w$ on $F'_{i,j}$ all lie in $V^+$ while their inneighbours lie in $V^-$.
For each $k=1,\dots,s$ we will denote the white subcluster of a cluster 
$U_k$ by $U_{k,w}$. We use similar notation for $U^+ _k$ and $U^- _k$.

Consider any $UV\in E(R')$. Recall that $U$ and $V$ are viewed as clusters of size $m$ in $R'$,
but when considering $\mathcal F_i$ we are in fact considering subclusters of $U$ and $V$ of size $m'$.
When viewed as clusters in $R'$, $UV$ initially corresponded to an $\eps$-regular pair of density at
least $\gamma_4 d'$ in $H_4$. Thus when viewed as clusters in $\mathcal F_i$, $UV$
initially corresponded to a $2\eps$-regular pair of density at least $\gamma _4 d'/2$ in $H_4$.
Moreover, initially the edges from $U_w$ to $V_w$ in $H_4$ induce a $3\eps$-regular pair of density
at least $\gamma _4 d'/3$. However, we have removed all the edges lying in the $\ell$ Hamilton cycles
$C_{i',j'}$ which we have defined already. Property~(iii) implies that
we have removed at most $3L^2 \ell\leq 3L^2n$ edges from $H_4$. Thus we have the following property:
\begin{itemize}
\item[$(\delta _1)$] Given any $UV \in E(R')$, let $S\subseteq U_w$, $T\subseteq V_w$ be such that
$|S|,|T|\geq 3\eps m'$. Then $e_{H_4} (S,T) \geq \gamma _4  d'|S||T|/4$.
\end{itemize}
When constructing $C_{i,j}$ we will remove at most $3L^2$ more edges from $H_4$.
But since $(\delta _1)$ is far from being tight, it will hold throughout the argument below.
Similarly, the initial definition of $H_{5,i}$ (c.f.~Section~\ref{randomsplit})
and~(iv) together imply the following property:
\begin{itemize}
\item[$(\delta _2)$] Consider any edge $V V^+\in E(\mathcal F_i)$. Let $S \subseteq V$ and
$T \subseteq V^+$ be such that $|S|,|T|\geq \sqrt{\eps} m'$.
Then $e_{H_{5,i}} (S,T) \geq \gamma _5 \beta |S||T|/2$.
\end{itemize}
We now construct $C_{i,j}$ from $F'_{i,j}$. Condition~$(\beta _6)$ in Section~\ref{4.6} implies that,
for each $k=1,\dots,s$, every vertex in $V_G(D'_k)$ lies on the same cycle, $C'_k$ say, in $F'_{i,j}$.
Let $x_1 \in U_{1,w}$ be such that $x_1$ has at least $\gamma _4 d'|U^+_{2,w}|/4\ge \gamma _4 d' m'/5$
outneighbours in $H_4$ which lie in $U_{2,w} ^+$. By~$(\delta _1)$ all but at most $3\eps m'$ vertices
in $U_{1,w}$ have this property. Note that the outneighbour in $F'_{i,j}$ of any such vertex lies in
$U^+ _1$. However, by~$(\delta _2)$ all but at most $\sqrt{\eps }m'$ vertices in $U^+ _1$
have at least $\gamma _5 \beta |U_{1,w}|/2\ge \gamma _5 \beta m'/3$ inneighbours in $H_{5,i}$ which lie in $U_{1,w}$.
Thus we can choose $x_1$ with the additional property that its outneighbour $y_1\in U^+_1$ in $F'_{i,j}$ has
at least $\gamma _5 \beta m'/3$ inneighbours in $H_{5,i}$ which lie in~$U_{1,w}$.

Let $P$ denote the directed path $C'_1-x_1 y_1$ from~$y_1$ to $x_1$. We now have
two cases to consider.

\medskip \noindent {\bf{Case 1.}} $C'_1 \not = C'_2$. 

\smallskip \noindent
Note that $x_1$ has at least $\gamma _4 d' m'/5-c'm'\gg \gamma _4 d' m'/6$ outneighbours
$y'_2\in U_{2,w} ^+$ in $ H_4$ such that the inneighbour of $y'_2$ in~$F'_{i,j}$ lies in $U_{2,w}$.
However, by~$(\delta_1)$ all but at most $3\eps m'$ vertices in  $U_{2,w}$ have 
at least $\gamma _4 d' m'/5$ outneighbours in $H_4$ which lie in $U^+ _{3,w}$.
Thus we can choose an outneighbour $y'_2\in U_{2,w} ^+$ of $x_1$ in $H_4$ such that
the inneighbour $x'_2$ of $y'_2$ in~$F'_{i,j}$ lies in $U_{2,w}$ and $x'_2$ has at least
$\gamma _4 d' m'/5$ outneighbours in $H_4$ which lie in $U^+ _{3,w}$. We extend $P$
by replacing it with $(P \cup C'_2 \cup \{ x_1 y'_2\}) \backslash \{x'_2 y'_2\}$.

\medskip \noindent {\bf{Case 2.}} $C'_1  = C'_2$.

\smallskip \noindent In this case the vertices in $V_G(D'_2)$ already lie on $P$.
We will use the following claim to modify $P$.

\begin{claim}\label{rotateclaim}
There is a vertex $y_2 \in U^+ _{2,w}$ such that:
\begin{itemize}
\item $x_1 y_2 \in E(H_4)$.
\item The predecessor $x_2$ of $y_2$ on $P$ lies in $U_{2,w}$.
\item There is an edge $x_2 y'_2$ in $H_{5,i}$ such that $y'_2 \in U^+ _{2,w}$ and $y_2$ precedes $y'_2$ on $P$
(but need not be its immediate predecessor).
\item The predecessor $x'_2$ of $y'_2$ on $P$ lies in $U_{2,w}$.
\item $x'_2$ has at least $\gamma _4 d' m'/5$ outneighbours in $H_4$ which lie in $U^+ _{3,w}$.
\end{itemize}
\end{claim}
\begin{figure}[htb!] 
\begin{center}\footnotesize
\psfrag{1}[][]{\normalsize $y_1$}
\psfrag{2}[][]{\normalsize $x_{2}$}
\psfrag{3}[][]{\normalsize $y_{2}$}
\psfrag{4}[][]{\normalsize $x' _{2}$}
\psfrag{5}[][]{\normalsize $y' _{2}$}
\psfrag{6}[][]{\normalsize $x_{1}$}
\includegraphics[width=0.65\columnwidth]{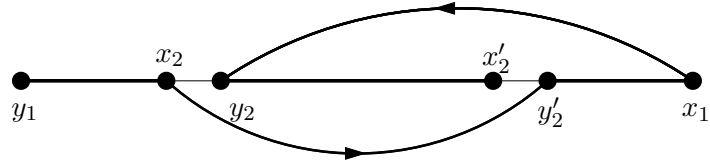}  
\caption{The modified path $P$ in Case~2} %
\end{center}
\end{figure}
\proof
Since $x_1 $ has at least $\gamma _4 d' m'/5$ outneighbours in $H_4$ which lie in $U^+ _{2,w}$,
at least $\gamma _4 d'm' /5- c' m'-3\eps m' \geq \gamma _4 d'm'/6$
of these outneighbours $y$ are such that the predecessor $x$ of $y$ on $P$ lies in $U_{2,w}$
and at least $\gamma _4 d'm' /5$ outneighbours of $x$ in $H_4$ lie in $U^+ _{3,w}$.
This follows since all such vertices $y$ have their predecessor on $P$ lying in $U_2$ (since $y \in U^+ _{2,w}$),
since $|U_{2,w}| \geq (1- c' )m'$ and since by~$(\delta_1)$ all but at most $3\eps m'$ vertices
in $U_{2,w}$ have at least $\gamma _4 d'm' /5$ outneighbours in $U^+ _{3,w}$.
Let $Y_2$ denote the set of all such vertices
$y$, and let $X_2$ denote the set of all such vertices $x$. So $|X_2|=|Y_2|\geq \gamma _4 d' m' /6$,
$X_2 \subseteq U_{2,w}$, $Y_2 \subseteq U^+ _{2,w} \cap N^+ _{H_4} (x_1)$.
Let $X^* _2 $ denote the set of the first $\gamma _4 d' m' /12$ vertices in $X_2$ on $P$ and
$Y^* _2$ the set of the last $\gamma _4 d' m'/12$ vertices in $Y_2$ on $P$.
Then~$(\delta_2)$ implies the existence of an edge $x_2y'_2$ from $X^*_2$ to $Y^*_2$ in~$H_{5,i}$.
Then the successor $y_2$ of $x_2$ on $P$ satisfies the claim.
\endproof

Let $x_2 , y_2, x'_2$ and $y'_2$ be as in Claim~\ref{rotateclaim}. We modify $P$ by replacing $P$ with
$$(P\cup\{ x_1 y_2, x_2 y'_2\}) \backslash \{ x_2 y_2 ,x'_2 y'_2\}$$
(see Figure~2).%

In either of the above cases we obtain a path $P$ from $y_1$ to some vertex $x'_2 \in U_{2,w}$ which has at
least $\gamma _4 d' m'/5$ outneighbours in $H_4$ lying in~$U^+ _{3,w}$.
We can repeat the above process: If $C'_3 \not = C'_1,C'_2$ then we extend $P$ as in Case~1.
If $C'_3 =C'_1$ or $C'_3 =C'_2$ then we modify $P$ as in Case~2. In both cases we obtain a
new path $P$ which starts in $y_1$ and ends in some $x'_3 \in U_{3,w}$ that has at least
$\gamma _4 d' m'/5$ outneighbours in $H_4$ lying in $U^+ _{4,w}$. 
We can continue this process, for each $C'_k$ in turn, until we obtain a path $P$ which
contains all the vertices in $C'_1, \dots , C'_s$ (and thus all the vertices in $G$),
starts in $y_1$ and ends in some $x'_s \in U_{s,w}$ having at least $\gamma _4 d' m'/5$ outneighbours
in $H_4$ which lie in $U^+ _{1,w}$.

\begin{claim}\label{rotateclaim2}
There is a vertex $y'_1 \in U^+ _{1}\setminus\{y_1\}$ such that:
\begin{itemize}
\item $x'_s y'_1 \in E(H_4)$.
\item The predecessor $x'_1$ of $y'_1$ on $P$ lies in $U_{1,w}$.
\item There is an edge $x'_1 y''_1$ in $H_{5,i}$ such that $y''_1 \in U^+ _{1,w}$ and $y'_1$ precedes
$y''_1$ on $P$.
\item The predecessor $x''_1$ of $y''_1$ on $P$ lies in $U_{1,w}$.
\item $x''_1$ has at least $\gamma _5\beta m'/3$ outneighbours in $H_{5,i}$ which lie in $U^+ _{1,w}$.
\end{itemize}
\end{claim}
\proof
The proof is almost identical to that of Claim~\ref{rotateclaim} except that we apply $(\delta_2)$ to ensure that
$x''_1$ has at least $\gamma _5\beta m'/3$ outneighbours in $H_{5,i}$ which lie in $U^+ _{1,w}$.
\endproof
Let $x'_1,y'_1, x''_1$ and $y''_1$ be as in Claim~\ref{rotateclaim2}. We modify $P$ by replacing it with
the path
$$(P\cup \{x'_sy'_1, x'_1 y''_1\} )\backslash \{x'_1 y'_1 , x''_1 y''_1\}$$
from $y_1$ to $x''_1$. So $P$ is a Hamilton path in~$G$ which is 
edge-disjoint from the $\ell$ Hamilton cycles $C_{i',j'}$ already defined. In each of
the $s$ steps in our construction of $P$ we have added at most one edge from
each of $H_4$ and $H_{5,i}$.  So by~(\ref{s}) $P$ contains at most $3L^2$ edges from $H_4$
and at most $3L^2$ edges from~$H_{5,i}$. All other edges of $P$ lie in $F'_{i,j}$. Recall that $y_1$
has at least $\gamma _5 \beta m'/3$ inneighbours in $H_{5,i}$ which lie in $U_{1,w}$ and $x''_1$ has 
at least $\gamma _5\beta m'/3$ outneighbours in $H_{5,i}$ which lie in $U^+ _{1,w}$.  Thus we can apply
Lemma~\ref{rotationlemma} to $P\cup H_{5,i}$ with $U^+_{1}$ playing the role of $V$ and
$U_{1}$ playing the role of $U$ to obtain a Hamilton cycle $C_{i,j}$ in $G$ where
$|E(C_{i,j})\backslash E(P)|\leq 5$. By construction, $C_{i,j}$ satisfies (i)--(iv). Thus we can
indeed find $(1/2-\eta _1)n $ Hamilton cycles in $G$, as desired.

\section{Almost decomposing oriented regular graphs with large semidegree} \label{38}

In this section, we describe how Theorem~\ref{main} can be extended to 
`almost regular' oriented graphs whose minimum semidegree is larger than $3n/8$.
More precisely,
we say that an oriented graph $G$ on $n$ vertices is \emph{$(\alpha  \pm \eta )n$-regular} if
$\delta^0(G)\ge (\alpha  - \eta) n$ and $\Delta^0(G) \le (\alpha  + \eta) n$.

\begin{thm}\label{main38} For every $ \gamma >0$ there exist $n_0=n_0 (\gamma )$ and
$\eta =\eta(\gamma)>0$ such that the following holds. Suppose
that $G$ is an $(\alpha \pm \eta)n$-regular oriented graph on $n\geq n_0$ vertices where 
$3/8+ \gamma \leq \alpha <1/2 $. Then $G$ contains at least $(\alpha -\gamma)n$
edge-disjoint Hamilton cycles.
\end{thm} 
Theorem~\ref{main38} is best possible in the sense that there
are almost regular oriented graphs whose semidegrees are all close to
$3n/8$ but which do not contain a Hamilton cycle. 
These were first found by H\"aggkvist~\cite{HaggkvistHamilton}.
However, we believe that if one requires $G$ to be completely regular, 
then one can actually obtain a Hamilton decomposition of $G$.
Note this would be a significant generalization of Kelly's conjecture.
\begin{conj}
For every $ \gamma >0$ there exists $n_0=n_0 (\gamma )$ such that
for all $n \ge n_0$ and all $r\ge (3/8+\gamma)n$ each    
$r$-regular oriented graph on $n$ vertices has a decomposition into Hamilton cycles.
\end{conj}
At present we do not even have any examples to rule out the possibility that
one can reduce the constant $3/8$ in the above conjecture:
\begin{question}
Is there a constant $c<3/8$ such that for every sufficiently large $n$ every $cn$-regular
oriented graph $G$ on $n$ vertices has a Hamilton decomposition or at least a set of edge-disjoint
Hamilton cycles covering almost all edges of $G$?
\end{question}
It is clear that we cannot take $c <1/4$ since 
there are non-Hamiltonian $k$-regular oriented graphs on $n$ vertices with 
$k=n/4-1/2$ (consider a union of 2 regular tournaments).

\medskip

\noindent
{\bf Sketch proof of Theorem~\ref{main38}.}
The proof of Theorem~\ref{main38} is similar to that of Theorem~\ref{main}.  
A detailed proof of Theorem~\ref{main38} can be found in~\cite{treglownphd}. 
The main use of the assumption of high minimum semidegree in our proof of 
Theorem~\ref{main} was that for any pair $A$, $B$ of large sets of vertices, we could assume the existence of many 
edges between $A$ and $B$ (see Lemma~\ref{keevashmult}). This enabled us to prove the existence of 
very short paths, shifted walks and skeleton walks between arbitrary pairs of vertices.
Lemma~\ref{keevashmult} does not hold under the weaker degree conditions of Theorem~\ref{main38}.
However, (e.g.~by Lemma~4.1 in~\cite{kelly}) these degree conditions are
strong enough to imply the following `expansion property':
for any set $S$ of vertices, we have that $|N^+_G(S)| \ge |S|+ \gamma n/2$
(provided $|S|$ is not too close to $n$).
Lemma~3.2 in~\cite{kelly} implies that this expansion property is also inherited by the reduced graph.
So in the proof of Lemma~\ref{multifactor1}, this expansion property can be used to find paths of
length $O(1/\gamma)$ which join up given pairs of vertices. Similarly, in
Lemma~\ref{shiftedwalk} we find closed shifted walks so that each cycle $C$ in $F$
is traversed $O(1/\gamma)$ times instead of just $3$ times (such a result is proved explicitly in Corollary~4.3
of~\cite{kelly}). Finally, in the proof of Claim~\ref{shadow} we now find shadow skeleton walks 
whose length is $O(1/\gamma)$ instead of 5.
In each of these cases, the increase in length does not affect the remainder of the proof.
\endproof

\section{Acknowledgement}
We would like to thank Demetres Christofides for helpful discussions.

\medskip 

{\footnotesize \obeylines \parindent=0pt

Daniela K\"{u}hn, Deryk Osthus \& Andrew Treglown
School of Mathematics
University of Birmingham
Edgbaston
Birmingham
B15 2TT
UK
}
\begin{flushleft}
{\it{E-mail addresses}:
\tt{\{kuehn,osthus,treglowa\}@maths.bham.ac.uk}}
\end{flushleft}

\end{document}